\documentclass[letterpaper,11pt]{amsart}
\usepackage{graphicx, amscd, amssymb}
\newtheorem{lem}{Lemma}[section]
\newtheorem{thm}[lem]{Theorem}
\newtheorem{pro}[lem]{Proposition}
\newtheorem{cor}[lem]{Corollary}

\newtheorem{con}[lem]{Conjecture}

\newcommand{\ZZ}{{\mathbb{Z}}}

\renewcommand{\P}{{\mathcal{P}}}

\newcommand{\R}{{\mathcal{R}}}

\newcommand{\M}{{\mathcal{M}}}
\newcommand{\Cat}{{\textrm{Cat}}}

\begin{document}

\title{The Cyclic Sieving Phenomenon for Faces of Generalized Cluster Complexes}

\author{Sen-Peng Eu}
\address{Department of
Applied Mathematics, National University of Kaohsiung, Kaohsiung
811, Taiwan, ROC} \email{speu@nuk.edu.tw}

\author{Tung-Shan Fu}
\address{Mathematics
Faculty, National Pingtung Institute of Commerce, Pingtung 900,
Taiwan, ROC} \email{tsfu@npic.edu.tw}

\thanks{The authors
are partially supported by National Science Council, Taiwan under
grants NSC 95-2115-M-390-006-MY3 (S.-P. Eu) and NSC
95-2115-M-215-001-MY2 (T.-S. Fu).}

\maketitle

\begin{abstract}
The notion of cyclic sieving phenomenon is introduced by Reiner,
Stanton, and White as a generalization of Stembridge's $q=-1$
phenomenon. The generalized cluster complexes associated to root
systems are given by Fomin and Reading as a generalization of the
cluster complexes found by Fomin and Zelevinsky. In this paper, the
faces of various dimensions of the generalized cluster complexes in
type $A_n$, $B_n$, $D_n$, and $I_2(a)$ are shown to exhibit the
cyclic sieving phenomenon under a cyclic group action. For the
cluster complexes of exceptional type $E_6$, $E_7$, $E_8$, $F_4$,
$H_3$, and $H_4$, a verification for such a phenomenon on their
maximal faces is given.
\end{abstract}

\section{Introduction}

 In \cite{RSW}, Reiner, Stanton, and White introduced the
notion of cyclic sieving phenomenon as a generalization of
Stembridge's $q=-1$ phenomenon for generating functions of a set of
combinatorial structures with a cyclic group action. Namely, a
triple $(X,X(q),C)$ consisting of a finite set $X$, a polynomial
$X(q)\in \ZZ[q]$ with the property that $X(1)=|X|$, and a cyclic
group $C$ that acts on $X$ is said to exhibit the {\em cyclic
sieving phenomenon} if for every $c\in C$,
\begin{equation} \label{eqn:Definition}
 [X(q)]_{q=\omega}=|\{x\in X: c(x)=x\}|,
\end{equation} where the complex number $\omega$ is
a root of unity of the same multiplicative order as $c$.
Equivalently, if $X(q)$ is expanded as $X(q)\equiv
a_0+a_1x+\cdots+a_{n-1}x^{n-1}$ (mod $q^n-1$), where $n$ is the
order of $C$, then $a_k$ counts the number of orbits on $X$ under
$C$, the stabilizer-order of which divides $k$. In particular, $a_0$
counts the total number of orbits, $a_1$ counts the number of free
orbits, and $a_2-a_1$ is the number of orbits that have a stabilizer
of order 2. This paper is motivated by the following concrete
example. Here we use the notation
${{n}\brack{i}}_q:=\frac{[n]!_q}{[i]!_q[n-i]!_q}$, where
$[n]!_q=[1]_q[2]_q\cdots[n]_q$ and $[i]_q=1+q+\cdots+q^{i-1}$.

\medskip
\begin{thm} \label{thm:motivation} {\rm (\cite[Theorem 7.1]{RSW})}
Let $X$ be the set of dissections of a regular $(n+2)$-gon using $k$
noncrossing diagonals $(0\le k\le n-1)$. Let
\begin{equation} \label{eqn:X(q)}
X(q):=\frac{1}{[k+1]_q}{{n+k+1}\brack{k}}_q{{n-1}\brack{k}}_q.
\end{equation}
Let the cyclic group $C$ of order $n+2$ act on $X$ by cyclic
rotation of the polygon. Then $(X,X(q),C)$ exhibits the cyclic
sieving phenomenon.
\end{thm}

Note that $X(1)=|X|$ is the well-known Kirkman-Cayley number. In
\cite{FZ}, Fomin and Zelevinsky  introduced a simplicial complex
$\Delta(\Phi)$, called the {\em cluster complex}, associated to a
root system $\Phi$, which can be realized by a combinatorial
structure constructed in terms of polygon-dissections. In fact,
Theorem \ref{thm:motivation} proves the cyclic sieving phenomenon
for the $k$-faces of the cluster complex $\Delta(\Phi)$ in type
$A_{n-1}$, under a cyclic group generated by a deformation $\Gamma$
(defined in Section 2) of Coxeter element of $\Phi$. This connection
will be explained in next section.

From \cite[Theorem 1.9]{FZ},  the number of {\em facets} (i.e.,
maximal faces) of $\Delta(\Phi)$ can be expressed uniformly as
\begin{equation} \label{eqn:facet1}
\Cat(\Phi):=\prod_{i=1}^n \frac{h+e_i+1}{e_i+1},
\end{equation}
known as the {\em generalized Catalan numbers} $\Cat(\Phi)$, where
$h$ is the Coxeter number and $e_1,\dots,e_n$ are the exponents of
$\Phi$. In particular, $\Cat(A_{n-1})=\frac{1}{n+1}{{2n}\choose{n}}$
is the $n$-th Catalan number. The cyclic group generated by $\Gamma$
is of order $h+2$. Along with the $q$-analogue of (\ref{eqn:facet1})
defined by
\begin{equation} \label{eqn:q-facet1} \Cat(\Phi,q):=\prod_{i=1}^n
\frac{[h+e_i+1]_q}{[e_i+1]_q},
\end{equation}
one of our main results is to prove the following theorem.

\medskip
\begin{thm} \label{thm:allCSP-facets} Let $X$ be the set of
facets of the cluster complex $\Delta(\Phi)$. Let
$X(q)=\Cat(\Phi,q)$ be defined in (\ref{eqn:q-facet1}). Let the
cyclic group $C$ of order $h+2$ generated by $\Gamma$ act on $X$.
Then $(X,X(q),C)$ exhibits the cyclic sieving phenomenon.
\end{thm}

\medskip
Moreover, in \cite{FR}, Fomin and Reading defined the {\em
generalized cluster complexes}
 $\Delta^s(\Phi)$ associated to a root
system $\Phi$ and a positive integer $s$, which specializes at $s=1$
to the cluster complexes $\Delta(\Phi)$. The purpose of this paper
is to study the cyclic sieving phenomenon for the faces of
$\Delta^s(\Phi)$, along with a $q$-analogue $X(q)$ of their face
numbers, under a cyclic group action. Making use of Fomin and
Reading's results, we prove the cyclic sieving phenomenon by a
combinatorial approach for $\Delta^s(\Phi)$ in type $A_n$, $B_n$,
$D_n$, and $I_2(a)$. For $\Phi$ of exceptional type $E_6$, $E_7$,
$E_8$, $F_4$, $H_3$, and $H_4$, although a systematic method is not
available, we verify such a phenomenon by computer for the facets of
$\Delta^s(\Phi)$ when $s=1$.

This paper is organized as follows. We review backgrounds of cluster
complexes and generalized cluster complexes in Section 2. As the
main results of this paper, the cyclic sieving phenomenon for the
generalized cluster complexes in type $A_n$, $B_n$, $D_n$, and
$I_2(a)$ are given in Sections 3, 4, 5, and 6, respectively. The
cases of exceptional type are given in Section 7. Finally, a
discussion regarding the polynomials $X(q)$ is given in Section 8.

\section{Backgrounds}

In this section, we review basic facts of cluster complexes from
\cite{FZ} and interprets Theorem \ref{thm:motivation} in terms of
cluster complexes (type $A_{n-1}$). Then we review the definition of
generalized cluster complexes from \cite{FR} and introduce the main
purpose of this paper. Most of this section follows the material in
\cite[Sections 2 and 3]{FR}

\subsection{Cluster complexes}
Let $\Phi$ be an irreducible root system of rank $n$. Let
$\Phi_{>0}$ denote the set of positive roots in $\Phi$, and let
$\Pi=\{\alpha_i:i\in I\}$ denote the set of simple roots in $\Phi$,
where $I=\{1,\dots,n\}$. Accordingly, $-\Pi=\{-\alpha_i:i\in I\}$ is
the set of negative simple roots. The set $S=\{s_i:i\in I\}$ of
reflections corresponding to simple roots $\alpha_i$ generates a
finite reflection group $W$ that naturally acts on $\Phi$. The pair
$(W,S)$ is a Coxeter system.

\begin{figure}[ht]
\begin{center}
\includegraphics[width=1.5in]{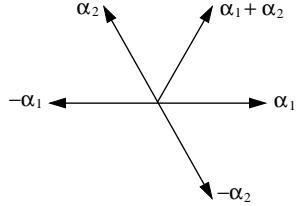}
\end{center}
\caption{\small The set $\Phi_{\ge -1}$ in type $A_2$.}
\label{fig:root-system}
\end{figure}

Let $I=I_+\cup I_-$ be a partition of $I$ such that each of sets
$I_+$ and $I_-$ is totally disconnected in the Coexter diagram. Let
$\Phi_{\ge -1}=\Phi_{>0}\cup (-\Pi)$. Define the involutions
$\tau_{\pm}:\Phi_{\ge -1}\rightarrow\Phi_{\ge -1}$ by
\[\tau_{\epsilon}(\alpha)=\left\{ \begin{array}{ll}
\alpha & \mbox{\rm if $\alpha=-\alpha_i$, for $i\in I_{-\epsilon}$},\\
\left(\prod_{i\in I_{\epsilon}} s_i\right)(\alpha) & \mbox{\rm
otherwise,}
       \end{array}
\right. \] for $\epsilon\in\{+,-\}$. The product
$\Gamma=\tau_-\tau_+$ generates a cyclic group
$\langle\Gamma\rangle$ that acts on $\Phi_{\ge -1}$. For example,
let $\Phi$ be the root system of type $A_2$, with $I=\{1,2\}$. The
set $\Phi_{\ge
-1}=\{\alpha_1,\alpha_2,\alpha_1+\alpha_2,-\alpha_1,-\alpha_2\}$ of
roots is shown in Figure \ref{fig:root-system}. Setting $I_+=\{1\}$
and $I_-=\{2\}$, one can check that $\Gamma$ acts on $\Phi_{\ge -1}$
by $\alpha_1\rightarrow
-\alpha_1\rightarrow\alpha_1+\alpha_2\rightarrow
-\alpha_2\rightarrow\alpha_2\rightarrow\alpha_1$.

By \cite[Section 3.1]{FZ}, the map $\Gamma$ induces a relation of
{\em compatibility} on $\Phi_{\ge -1}$ such that (i)
$\alpha,\beta\in\Phi_{\ge -1}$ are compatible if and only if
$\Gamma(\alpha)$ and $\Gamma(\beta)$ are compatible; (ii)
$-\alpha_i\in -\Pi$ and $\beta\in \Phi_{>0}$ are compatible if and
only if the simple root expansion of $\beta$ does not involve
$\alpha_i$. Following \cite[p. 983]{FZ}, the {\em cluster complex}
$\Delta(\Phi)$ is defined to be the simplicial complex whose faces
are subsets of roots in $\Phi_{\ge -1}$, which are pairwise
compatible. For $\Phi$ in type $A_n$, $B_n$, and $D_n$, the cluster
complex $\Delta(\Phi)$ can be realized by dissections of a regular
polygon such that $\langle\Gamma\rangle$ corresponds to a group
action on the dissections by rotation of the given polygon.

\begin{figure}[ht]
\begin{center}
\includegraphics[width=1.5in]{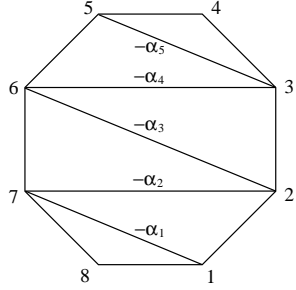}
\end{center}
\caption{\small The snake in type $A_5$.} \label{fig:snake}
\end{figure}

\begin{figure}[ht]
\begin{center}
\includegraphics[width=1.5in]{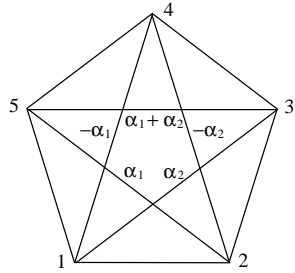}
\end{center}
\caption{\small A representation for the set $\Phi_{\ge -1}$ of
roots in type $A_2$.} \label{fig:realization}
\end{figure}

Specifically, consider the root system $\Phi$ of type $A_n$. Then
$\Phi_{\ge -1}$ consists of positive roots of the form
$\alpha_{ij}=\alpha_i+\alpha_{i+1}+\cdots+\alpha_j$, for $1\le i\le
j\le n$, and negative simple roots $-\alpha_i$, for $1\le i\le n$.
Let $P$ be a regular polygon with $n+3$ vertices labeled by
$\{1,2,\dots,n+3\}$ counterclockwise. The roots in $\Phi_{\ge -1}$
are identified with diagonals of $P$ as follows. For $1\le
i\le\frac{n+1}{2}$, the root $-\alpha_{2i-1}\in -\Pi$ is identified
with the diagonal connecting vertices $i$ and $n+3-i$. For $1\le
i\le \frac{n}{2}$, the root $-\alpha_{2i}\in -\Pi$ is identified
with the diagonal connecting vertices $i+1$ and $n+3-i$. These
diagonals form a `snake' of negative simple roots. Figure
\ref{fig:snake} shows this snake in type $A_5$. The positive roots
of $\Phi_{>0}$ are identified with the remaining diagonals as
follows. Each $\alpha_{ij}$ is identified with the unique diagonal
that intersects the diagonals
$-\alpha_i,-\alpha_{i+1},\dots,-\alpha_j$ and no other diagonals in
the snake. For example, Figure \ref{fig:realization} is a
realization of the set $\Phi_{\ge -1}$ in type $A_2$.

Under this bijection, every pair of compatible roots is carried to a
pair of noncrossing diagonals. Hence each $k$-face (i.e.,
$k$-element simplex) of the cluster complex $\Delta(A_n)$
corresponds to a dissection of $P$ using $k$ noncrossing diagonals.
Moreover, the map $\Gamma$ corresponds to a clockwise rotation of
$P$ that carries point $2$ to point $1$, etc. Therefore, Theorem
\ref{thm:motivation} can be interpreted in terms of the cluster
complex $\Delta(A_{n-1})$, i.e., let $X$ be the set of $k$-faces of
$\Delta(A_{n-1})$, let $X(q)$ be the polynomial defined by
(\ref{eqn:X(q)}). Then $(X,X(q),\langle\Gamma\rangle)$ exhibits the
cyclic sieving phenomenon.

\subsection{Generalized cluster complexes}
Let $s$ be a positive integer. For each $\alpha\in\Phi_{>0}$, let
$\alpha^1,\dots,\alpha^s$ denote the $s$ `colored' copies of
$\alpha$. Define
\[\Phi^s_{\ge -1}=\{\alpha^k:\alpha\in\Phi_{>0},1\le k\le
s\}\cup\{(-\alpha_i)^1:i\in I\}, \]
 i.e., $\Phi^s_{\ge -1}$ consists
of $s$ copies of the positive roots and one copy of the negative
simple roots. The relation of compatibility on $\Phi^s_{\ge -1}$ can
be defined by an $s$-analogue of $\Gamma$. For
$\alpha^k\in\Phi^s_{\ge -1}$, define
\[\Gamma_s(\alpha^k)=\left\{ \begin{array}{ll}
\alpha^{k+1} & \mbox{\rm if $\alpha\in\Phi_{>0}$ and $k<s$},\\
(\Gamma(\alpha))^1 & \mbox{\rm otherwise.}
       \end{array}
\right. \] By \cite[Theorem 3.4]{FR}, the map $\Gamma_s$ induces a
relation of compatibility on $\Phi^s_{\ge -1}$ such that (i)
$\alpha^k$ and $\beta^l$ are compatible if and only if
$\Gamma_s(\alpha^k)$ and $\Gamma_s(\beta^l)$ are compatible; (ii)
$(-\alpha_i)^1$ and $\beta^l$ are compatible if and only if the
simple root expansion of $\beta^l$ does not involve $\alpha_i$. The
{\em generalized cluster complex} $\Delta^s(\Phi)$ associated to a
root system $\Phi$ is defined to be the simplicial complex whose
faces are subsets of roots in $\Phi^s_{\ge -1}$, which are pairwise
compatible.

A feasible $q$-analogue $X(q)$ of face numbers plays an essential
role in the cyclic sieving phenomenon. Fomin and Reading derived a
unified formula for the face numbers of various dimensions of the
complexes $\Delta^s(\Phi)$ in terms of Coxeter numbers and exponents
(see \cite[Theorem 8.5]{FR}). Let $f_k(\Phi,s)$ denote the $k$-{\em
face number} of the complex $\Delta^s(\Phi)$. For $\Phi$ of type
$A_n$, $B_n$, $D_n$, and $I_2(a)$, the $k$-face numbers can be
expressed explicitly as follows.

\medskip
\begin{thm} \label{thm:formula} For $\Phi$ of type $A_{n}$, $B_n$, $D_n$, and
$I_2(a)$, the face numbers of the generalized cluster complex
$\Delta^s(\Phi)$ are given by
\begin{enumerate}
\item
${\displaystyle
f_k(A_{n},s)=\frac{1}{k+1}{{s(n+1)+k+1}\choose{k}}{{n}\choose{k}}}$,
\\

\item ${\displaystyle f_k(B_n,s)={{sn+k}\choose{k}}{{n}\choose{k}}}$, \\

\item
${\displaystyle f_k(D_n,s)={{s(n-1)+k}\choose{k}}{{n}\choose{k}}+{{s(n-1)+k-1}\choose{k}}{{n-2}\choose{k-2}}}$,\\

\item $f_1(I_2(a),s)=sa+2$, and $f_2(I_2(a),s)={\displaystyle \frac{(sa+2)(s+1)}{2}}$.

\end{enumerate}

\end{thm}

\medskip
Case (i) of Theorem \ref{thm:formula} is due to J.~H.~Przytycki and
A.~S.~Sikora \cite{PS} and case (ii) is due to E.~Tzanaki
\cite{Tz-JCTA}. We remark that an obvious $q$-analogue of Fomin and
Reading's unified formula is not a feasible option for the
polynomial $X(q)$. (Some cases are even not polynomials with
nonnegative integer coefficients). The $q$-analogue $X(q)$ that
serves our purpose are derived case by case from Theorem
\ref{thm:formula}. (See Theorems \ref{thm:CSP}, \ref{thm:CSP-B},
\ref{thm:CSP-D}, and \ref{thm:CSP-I2(a)}).

\medskip
For the special case $k=n$, the number of maximal faces of the
complex $\Delta^s(\Phi)$ can be expressed uniformly as
\begin{equation} \label{eqn:facet}
\Cat^{(s)}(\Phi):=\prod_{i=1}^n \frac{sh+e_i+1}{e_i+1},
\end{equation}
where the {\em Coxeter number} $h$ and the {\em exponents}
$e_1,\dots,e_n$ of $\Phi$ are listed in Figure \ref{fig:table}. Note
that the expression (\ref{eqn:facet}) specializes at $s=1$ to the
generalized Catalan numbers (\ref{eqn:facet1}).

\begin{figure}[ht]
\begin{tabular}{ccccl}
  \hline
  % after \\: \hline or \cline{col1-col2} \cline{col3-col4} ...
  $\Phi$ & & $h$ & & $e_1,\dots,e_n$ \\
  \hline
  $A_n$ & & $n+1$ & & $1,2,\dots,n$ \\
  $B_n$ & & $2n$  & & $1,3,\dots,2n-1$ \\
  $D_n$ & & $2(n-1)$ & & $1,3,\dots,2n-3,n-1$ \\
  $E_6$ & & $12$  & & $1,4,5,7,8,11$ \\
  $E_7$ & & $18$  & & $1,5,7,9,11,13,17$ \\
  $E_8$ & & $32$  & & $1,7,11,13,17,19,23,29$ \\
  $F_4$ & & $12$  & & $1,5,7,11$ \\
  $H_3$ & & $10$  & & $1,5,9$ \\
  $H_4$ & & $30$  & & $1,11,19,29$ \\
  $I_2(a)$ & & $a$ & & $1,a-1$ \\
  \hline
\end{tabular}
\caption{\small The Coxeter number and exponents of $\Phi$.}
\label{fig:table}
\end{figure}

In this paper, we aim to prove the cyclic sieving phenomenon for the
generalized cluster complexes $\Delta^s(\Phi)$ in the framework that
$X$ is the set of $k$-faces of $\Delta^s(\Phi)$, $C$ is the cyclic
group of order $sh+2$ generated by $\Gamma_s$, and $X(q)$ is a
$q$-analogue of the $k$-face numbers.  For $\Phi$ of type $A_n$,
$B_n$, and $D_n$, our results rely on Fomin and Reading's
realization constructed in terms of polygon-dissections
\cite[Section 5]{FR}. Under this realization, the cyclic group $C$
corresponds to rotation of the given polygon. In type $I_2(a)$, we
make use of the graph-representation of the complex
$\Delta^s(I_2(a))$ given in \cite[Example 4.4]{FR}. For $\Phi$ of
exceptional type $E_6$, $E_7$, $E_8$, $F_4$, $H_3$, and $H_4$, a
complete verification of such a phenomenon only for the maximal
faces of $\Delta^s(\Phi)$ and only when $s=1$ is given.

\section{The cyclic sieving phenomenon for $\Delta^s(A_{n-1})$}

\begin{figure}[ht]
\begin{center}
\includegraphics[width=2in]{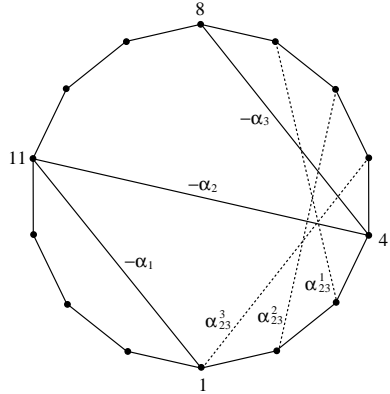}
\end{center}
\caption{\small The $s$-snake in type $A_3$.}
\label{fig:s-realization}
\end{figure}

Let $P$ be a regular polygon with $sn+2$ vertices labeled by
$\{1,2,\dots,sn+2\}$ counterclockwise. Consider the set of
dissections of $P$ into $(sj+2)$-gons $(1\le j\le n-1$) by
noncrossing diagonals. Such dissections are called $s$-{\em
divisible}. For convenience, a diagonal in an $s$-divisible
dissection is called $s$-{\em divisible}. Consider a root system
$\Phi$ of type $A_{n-1}$. Following \cite[Section 5.1]{FR}, the
roots of $\Phi^s_{\ge -1}$ can be identified with the $s$-divisible
diagonals of $P$ as follows. For $1\le i\le\frac{n}{2}$, the root
$-\alpha_{2i-1}$ is identified with the diagonal connecting points
$s(i-1)+1$ and $s(n-i)+2$. For $1\le i\le \frac{n-1}{2}$, the root
$-\alpha_{2i}$ is identified with the diagonal connecting points
$si+1$ and $s(n-i)+2$. These $n-1$ diagonals form an $s$-snake of
negative simple roots. For each positive root
$\alpha_{ij}=\alpha_i+\cdots+\alpha_j$ ($1\le i\le j\le n-1$), there
are exactly $s$ diagonals, which are $s$-divisible, intersecting the
diagonals $-\alpha_i,\dots,-\alpha_j$ and no other diagonals in the
$s$-snake. This collection of diagonals is of the form $D$,
$\Gamma_s^1(D),\dots,\Gamma_s^{s-1}(D)$ for some diagonal $D$. For
$1\le k\le s$, we identify $\alpha_{ij}^k$ with $\Gamma_s^{k-1}(D)$.
Figure \ref{fig:s-realization} shows the $s$-snake for $s=3$ and
$n=4$, along with the diagonals identified with the colored roots
$\alpha_{23}^1$, $\alpha_{23}^2$, and $\alpha_{23}^3$. Under this
bijection, the $k$-faces of the complex $\Delta^s(A_{n-1})$
correspond to the $s$-divisible dissections of $P$ using $k$
noncrossing diagonals, and $\Gamma_s$ corresponds to clockwise
rotation of $P$ carrying point $2$ to point $1$, etc.

A feasible polynomial for $X(q)$ is the natural $q$-analogue of
Theorem \ref{thm:formula}(i). Define
\begin{equation} \label{eqn:G(q)}
G(s,n,k;q)=\frac{1}{[k+1]_q}{{sn+k+1}\brack{k}}_q{{n-1}\brack{k}}_q,
\end{equation}
for $0\le k\le n-1$. Note that $G(s,n,k;1)=f_k(A_{n-1},s)$. As a
generalization of \cite[Theorem 7.1]{RSW}, we prove that the
$k$-faces of $\Delta^s(A_{n-1})$ exhibit the cyclic sieving
phenomenon under the group action $\langle\Gamma_s\rangle$.

\medskip
\begin{thm} \label{thm:CSP} For positive integers $s$ and $n$, let $X$ be the set of
$s$-divisible dissections of an $(sn+2)$-gon using $k$ noncrossing
diagonals. Let the cyclic group $C$ of order $sn+2$ act on $X$ by
cyclic rotation of the polygon. Let $X(q):=G(s,n,k;q)$. Then
$(X,X(q),C)$ exhibits the cyclic sieving phenomenon.
\end{thm}

\begin{figure}[ht]
\begin{center}
\includegraphics[width=4.4in]{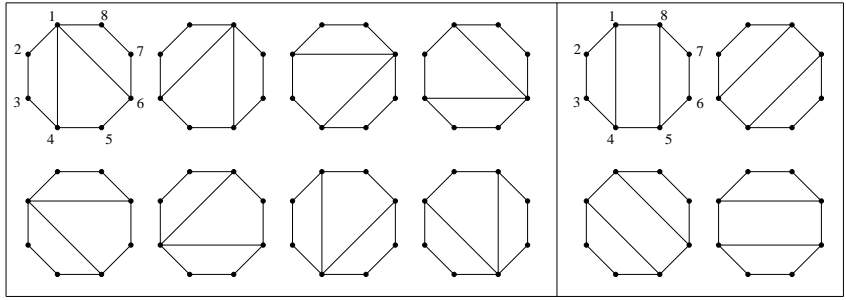}
\end{center}
\caption{\small The 2-divisible dissections of an octagon using 2
noncrossing diagonals.} \label{fig:octagon}
\end{figure}

\medskip For example, take $s=2$, $n=3$, and $k=2$. Then
$X(q)\equiv 2+q+2q^2+q^3+2q^4+q^5+2q^6+q^7$ (mod $q^8-1$). As shown
in Figure \ref{fig:octagon}, there are 12 $2$-divisible dissections
of an octagon using 2 noncrossing diagonals. These dissections are
partitioned into two orbits under a group action by cyclic rotation,
one of which is free and the other has a stabilizer of order 2.

We shall prove Theorem \ref{thm:CSP} by verifying condition
(\ref{eqn:Definition}) mentioned in the introduction. Recall that a
complex $\omega$ is a primitive $d$-th root of unity if $d$ is the
least integer such that $\omega^d=1$. First, we evaluate $X(q)$ at
primitive roots of unity (Proposition \ref{pro:evaluation}).  To do
this, we make use of the following facts and the $q$-Lucas theorem
(see \cite[Proposition 2.1]{GZ}).

\medskip
\begin{lem} \label{lem:Lemma} Let $m$, $m_1$, $m_2$, $k$, and $d$ be positive
integers. Let $\omega$ be a primitive $d$-th root of unity. Then
\begin{enumerate}
\item $[m]_{q=\omega}=0$ if and only if $d$ divides $m$ ($d\neq
  1$).

\item If $m_1\equiv m_2$ (mod $d$), then
 \[
  \lim_{q\rightarrow\omega} \frac{[m_1]_q}{[m_2]_q}=
 \left\{
\begin{array}{cl} \frac{m_1}{m_2} & \mbox{\rm if $m_1\equiv m_2\equiv 0$ (mod $d$)},\\
1 & \mbox{\rm if $m_1\equiv m_2\not\equiv 0$ (mod $d$).}
       \end{array}
\right.
\]

\item {\rm ($q$-Lucas Theorem)} If $m=ad+b$ and $k=rd+s$, where $0\le r,s\le d-1$,  then
\[{{m}\brack{k}}_{q=\omega}={{a}\choose{r}}{{b}\brack{s}}_{q=\omega}.
\]

\end{enumerate}

\end{lem}

\medskip
As a special case of Lemma \ref{lem:Lemma}(iii), the following
result is repeatedly used.

\medskip
\begin{cor}\label{cor:d|k} {\rm(\cite{BFW})} Let $m$, $k$, and $d$ be positive integers.
Let $\omega$ be a primitive $d$-th root of unity. For $d\ge 2$ a divisor of $m$,  then
\[{{m+k-1}\brack{k}}_{q=\omega}=\left\{ \begin{array}{cl}

\vspace{0.15cm}

{\displaystyle{{\frac{m+k}{d}-1}\choose{\frac{k}{d}}}}
& \mbox{if $d|k$,}\\

0 & \mbox{otherwise.}
\end{array}
 \right.
\]
\end{cor}

\medskip
\begin{pro} \label{pro:evaluation}
For $d\ge 2$ a divisor of $sn+2$, let $\omega$ be a  primitive
$d$-th root of unity. Then
\[
[G(s,n,k;q)]_{q=\omega}=\left\{ \begin{array}{cl}

\vspace{0.15cm}

{\displaystyle
{{\frac{sn+k+1}{2}}\choose{\frac{k+1}{2}}}{{\frac{n-2}{2}}\choose{\frac{k-1}{2}}}}
& \mbox{if $d=2$, $k$ odd, and $n$ even}\\

\vspace{0.15cm}

{\displaystyle {{\frac{sn+2+k}{d}-1}\choose{\frac{k}{d}}}
                {{\lfloor\frac{n-1}{d}\rfloor}\choose{\frac{k}{d}}}}
& \mbox{if $d\ge 2$ and $d|k$,} \\
0 & \mbox{otherwise.}
\end{array}
 \right.
\]
\end{pro}

\medskip
\begin{proof} Evaluating $G(s,n,k;q)$ at $q=\omega$, we have
\begin{eqnarray}
[G(s,n,k;q)]_{q=\omega} &=& \lim_{q\rightarrow\omega}
\frac{1}{[k+1]_{q}}{{sn+k+1}\brack{k}}_{q}{{n-1}\brack{k}}_{q} \nonumber\\
&=& \lim_{q\rightarrow\omega} \frac{([sn+k+1]_q\cdots
[sn+2]_q)([n-1]_q \cdots[n-k]_q)}{[k+1]!_q[k]!_q}.
\label{eqn:parity}
\end{eqnarray}
For $d=2$, making use of (i) and (ii) of Lemma \ref{lem:Lemma} and
examining parities of the factors in the numerator and denominator
in (\ref{eqn:parity}), we have
\begin{equation} \label{eqn:d<3}
[G(s,n,k;q)]_{q=\omega} =\left\{ \begin{array}{cl} \vspace{0.15cm}
{{\frac{sn+k+1}{2}}\choose{\frac{k+1}{2}}}{{\frac{n-2}{2}}\choose{\frac{k-1}{2}}}
& \mbox{if $k$ odd and $n$  even,} \\
{{\frac{sn+k}{2}}\choose{\frac{k}{2}}}{{\lfloor\frac{n-1}{2}\rfloor}\choose{\frac{k}{2}}}
& \mbox{if $k$  even,}\\
0 & \mbox{otherwise.}
\end{array}
\right.
\end{equation}

For $d\ge 3$, by Corollary \ref{cor:d|k}, we have
\[\lim_{q\rightarrow\omega}
{{sn+k+1}\brack{k}}_{q} =\left\{
\begin{array}{cl} \vspace{0.15cm}
{{\frac{sn+2+k}{d}-1}\choose{\frac{k}{d}}}
& \mbox{if $d|k$,}\\
0 & \mbox{otherwise.}
\end{array}
\right.
\]
For the rest of $[G(s,n,k;q)]_{q=\omega}$, we have
\begin{equation} \label{eqn:infinity}
\lim_{q\rightarrow\omega}
\frac{1}{[k+1]_q}{{n-1}\brack{k}}_{q}=\lim_{q\rightarrow\omega}\frac{[n-1]_q\cdots[n-k]_q}{[k+1]!_q}.
\end{equation}
Note that $d\nmid n$ since $d$ divides $sn+2$. For $0\le k\le n-1$,
one can check that there are at least as many factors $[t]_q$ in the
numerator as the denominator in (\ref{eqn:infinity}) such that $d$
divides $t$, and hence (\ref{eqn:infinity}) does not tend to
infinity. Moreover, for $d|k$, $\lim_{q\rightarrow\omega} [k+1]_q=1$
and
\[
\lim_{q\rightarrow\omega} \frac{1}{[k+1]_q}{{n-1}\brack{k}}_{q}=
{{\lfloor\frac{n-1}{d}\rfloor}\choose{\frac{k}{d}}}.
\]
It follows that, for $d\ge 3$, we have
\begin{equation} \label{eqn:d>=3}
[G(s,n,k;q)]_{q=\omega}=\left\{
\begin{array}{cl} \vspace{0.15cm}
{{\frac{sn+2+k}{d}-1}\choose{\frac{k}{d}}}{{\lfloor\frac{n-1}{d}\rfloor}\choose{\frac{k}{d}}}
& \mbox{if $d|k$,}\\
0 & \mbox{otherwise.}
\end{array}
\right.
\end{equation}
The assertion follows from (\ref{eqn:d<3}) and (\ref{eqn:d>=3}).
\end{proof}

\medskip Let $X$ be the set of $s$-divisible dissections of an
$(sn+2)$-gon $P$ using $k$ noncrossing diagonals. Let $C$ be the
cyclic group of order $sn+2$ acting on $X$ by cyclic rotation of
$P$. For $d\ge 2$ a divisor of $sn+2$, let $U(s,n,k,d)\subseteq X$
denote the set of dissections that are invariant under $d$-fold
rotation (i.e., a subgroup $C_d$ of order $d$ of $C$ generated by a
$\frac{2\pi}{d}$-rotation of $P$). In the following, we shall
enumerate $U(s,n,k,d)$ to complete the proof of Theorem
\ref{thm:CSP}.

For any centrally symmetric dissection, we observe that there is
either a diameter or an $(sj+2)$-gon in the center.  Hence if
$U(s,n,k,d)$ is nonempty, then either $d=2$ and $k$ is odd, or $d\ge
2$ and $d$ divides $k$. These two cases are treated in Propositions
\ref{pro:k-odd} and \ref{pro:bijection}, respectively.

\medskip
\begin{pro} \label{pro:k-odd} For positive integers $s$ and $n$, with $sn$ even and $k$ odd,
we have
\[ |U(s,n,k,2)|=
\left\{ \begin{array}{cl} \vspace{0.15cm}

{\displaystyle
{{\frac{sn+k+1}{2}}\choose{\frac{k+1}{2}}}{{\frac{n-2}{2}}\choose{\frac{k-1}{2}}}}
& \mbox{if $n$ even,}\\

0 & \mbox{otherwise.}
       \end{array}
\right.
\]
\end{pro}

\medskip
\begin{proof}
Note that a dissection with odd number of diagonals, which is
invariant under $2$-fold rotation, contains a unique diameter. This
diameter dissects the polygon into two $(sm+2)$-gons, where
$m=\frac{n}{2}$. Hence $n$ is even. Since there are $\frac{sn+2}{2}$
diameters and the dissection is completely determined by the
subdivision of the $(sm+2)$-gon on either side of the diameter using
$\frac{k-1}{2}$ diagonals, by (\ref{eqn:G(q)}) the number of such
dissections is
\[\frac{sn+2}{2}\cdot G\left(s,\frac{n}{2},\frac{k-1}{2};1\right)=
{{\frac{sn+k+1}{2}}\choose{\frac{k+1}{2}}}{{\frac{n-2}{2}}\choose{\frac{k-1}{2}}},
\]
as required.
\end{proof}

\medskip
For the latter case, the result relies on a bijection (Proposition
\ref{pro:bijection}), which is inspired by work of Tzanaki
\cite{Tz-JCTA}. In fact, for the special case $d=2$ and $n$ even,
the result has been obtained by Tzanaki in \cite[Corollary
3.2]{Tz-JCTA} by a bijection similar to the one given by Przytycki
and Sikora in \cite[Theorem 1]{PS}. We extend this method to
enumerate $d$-fold rotationally symmetric dissections for all $d\ge
2$.

\medskip
\begin{pro} \label{pro:bijection} For $d\ge 2$ a common divisor of $sn+2$ and $k$,
there is a bijection between the set $U(s,n,k,d)$ and the cartesian
product of the set of sequences $\{1\le a_1\le a_2\le\cdots\le
a_\frac{k}{d}\le\frac{sn+2}{d}\}$ and the set of sequences
$(\epsilon_1,\dots,\epsilon_m)\in\{0,1\}^m$ with exactly
$\frac{k}{d}$ entries equal to $1$, where
$m={\lfloor\frac{n-1}{d}\rfloor}$.
\end{pro}

\medskip
\begin{proof} Fixing $s$ and $d$, we define
$N_r:=\{n\in\ZZ: n-1\equiv r$ (mod $d$), and $d|(sn+2)\}$, for $0\le
r\le d-1$. Given an $N_r$ and an $n\in N_r$, let $b=\frac{sn+2}{d}$.
The vertices of an $(sn+2)$-gon are partitioned into $d$ {\em
sectors} of size $b$ by a vertex-labeling $\{1_1,2_1,\dots,b_1\}$,
$\{1_2,2_2,\dots,b_2\},\dots,\{1_d,2_d,\dots,b_d\}$ in the
counterclockwise order. For each $\pi\in U(s,n,k,d)$, there is no
diameter in $\pi$ since $d$ divides $k$. Every diagonal of $\pi$ is
oriented in such a way that if you travel along it in this direction
then the center of the polygon is on the left. These diagonals are
partitioned into $\frac{k}{d}$ orbits of size $d$ under the $d$-fold
rotation $C_d$ such that each orbit $R$ contains a unique diagonal
that leaves from the first sector $\{1_1,2_1,\dots,b_1\}$. We
associate $R$ with the label $i\in [b]$ if that diagonal leaves from
$i_1$.

Let $m=\frac{n-1-r}{d}$. Define
\begin{eqnarray*}
& & A(m,k)=\{(a_1,\dots,a_{\frac{k}{d}}):1\le a_1\le\cdots\le
a_{\frac{k}{d}}\le b\}, \\
& & B(m,k)=\{(\epsilon_1,\dots,\epsilon_m)\in \{0,1\}^m:\mbox{
exactly $\frac{k}{d}$ entries $\epsilon_j$ are equal to 1$\}$.}
\end{eqnarray*}
We shall establish a bijection $\Phi_m: U(s,n,k,d)\rightarrow
\{(\mu,\nu)|\mu\in A(m,k), \nu\in B(m,k)\}$
 by induction on $m$. Recall that $k\le n-1$ and $d$ divides $k$. If $m=1$ then
$\frac{k}{d}< 2$. For $k=0$, the map $\Phi_1$ sends the trivial
dissection to $(\emptyset, (\epsilon_1))$, where $\epsilon_1=0$, and
for $\frac{k}{d}=1$, $\Phi_1$ sends each dissection $\pi$ to
$((a_1),(\epsilon_1))$, where $a_1$ is the label of the unique
diagonal-orbit in $\pi$ and $\epsilon_1=1$.

For the inductive step, assume that we have established the
bijection $\Phi_{m-1}$, for $0\le \frac{k}{d}\le m-1$. To determine
$\Phi_{m}$, let $P$ be an $(sn+2)$-gon, and let  $\pi\in
U(s,n,k,d)$. The result is trivial for $k=0$. For $\frac{k}{d}\ge
1$, let $\mu=(a_1,\dots,a_{\frac{k}{d}})$ be the sequence of labels
associated to the diagonal-orbits in $\pi$, where $1\le
a_1\le\cdots\le a_{\frac{k}{d}}\le b$. Locate the first $a_j$ $(1\le
j\le \frac{k}{d}$) such that $a_j+s+1\le a_{j+1}$ (where
$a_{\frac{k}{d}+1}$ is understood to mean $a_1+b$). Such $a_j$
exists by pigeon hole principle. With abuse of notation, we identify
$a_i$ with the vertex labeled by $(a_i)_1$ and $a_i+k$ with the
vertex at distance $k$ from $a_i$ in the counterclockwise order. We
set $\epsilon_1=1$ if $a_j(a_j+s+1)$ is a diagonal in $\pi$, and
$\epsilon_1=0$ otherwise. Let $P'$ be the truncated $(sn+2-sd)$-gon
obtained from $P$ by removing the vertex-orbit of the set
$\{a_j+1,\dots,a_j+s\}$ under $C_d$, which is dissected by the
remaining $\frac{k}{d}-\epsilon_1$ orbits of diagonals, and let
$\pi'$ be the dissection determined by these orbits. Note that
$sn+2-sd\in N_r$ and that the labels associated to these
diagonal-orbits in $\pi'$ are the same, whether considered in $P$ or
$P'$. Therefore, by induction, $\Phi_{m-1}$ carries $\pi'$ to a pair
$(\mu',\nu')\in A(m-1,k)\times B(m-1,k)$, where
$\mu'=\{a_1,\dots,a_{\frac{k}{d}}\}$ if $\epsilon_1=0$ and
$\mu'=\{a_1,\dots,a_{j-1},a_{j+1},\dots,a_{\frac{k}{d}}\}$
otherwise, and $\nu'=(\epsilon_2,\dots,\epsilon_{m})$ with
$\frac{k}{d}-\epsilon_1$ entries $\epsilon_j$ equal to 1. We define
$\Phi_{m}(\pi)$ to be the pair
$((a_1,\dots,a_{\frac{k}{d}}),(\epsilon_1,\dots,\epsilon_{m}))$.

To define $\Phi_m^{-1}$, we proceed by induction. The case $m=1$ is
obvious. Assume that $\Phi_{m-1}^{-1}$ has been determined. Given
the sequences $\mu=(a_1,\dots,a_{\frac{k}{d}})$, where $1\le
a_1\le\cdots\le a_{\frac{k}{d}}\le b$, and
$\nu=(\epsilon_1,\dots,\epsilon_{m})$ with $\frac{k}{d}$ entries
equal to 1, we shall retrieve the dissection
$\Phi_{m}^{-1}((\mu,\nu))$ of an $(sn+2)$-gon $P$ from $(\mu,\nu)$
by recovering the set of diagonal-orbits. For $k=0$, the element
$(\emptyset,(0,\dots,0))$ is carried to the trivial dissection. For
$\frac{k}{d}\ge 1$, locate the first $a_j$ $(1\le j\le \frac{k}{d}$)
such that $a_j+s+1\le a_{j+1}$. Then there is an orbit $R$
containing the diagonal $a_j(a_j+s+1)$ if $\epsilon_1=1$ and there
is not otherwise. Let $P'$ be the $(sn+2-sd)$-gon the $i$-th sector
of which consists the vertices
$\{1_i,\dots,(a_j)_i,(a_j+s+1)_i,\dots,b_i\}$, for $1\le i\le d$.
Let $\mu'=(a_1,\dots,a_{\frac{k}{d}})$ if $\epsilon_1=0$ and
$\mu'=(a_1,\dots,a_{j-1},a_{j+1},\dots,a_{\frac{k}{d}})$ otherwise,
and let $\nu'=(\epsilon_2,\dots,\epsilon_{m})$ with
$\frac{k}{d}-\epsilon_1$ entries $\epsilon_j$ equal to 1. By
induction, the set $\R$ of diagonal-orbits  in $\P'$ is recovered
from $(\mu',\nu')$ by $\Phi_{m-1}^{-1}$. Therefore, we define
$\Phi_{m}^{-1}((\mu,\nu))$ to be the dissection determined by $\R$
if $\epsilon_1=0$ and by $\R\cup\{R\}$ if $\epsilon_1=1$.

Hence we prove the assertion for the $(sn+2)$-gons, where $n\in
N_r$. Since the argument works well for all $N_r$ ($0\le r\le d-1$),
the proof is completed.
\end{proof}

\begin{figure}[ht]
\begin{center}
\includegraphics[width=4in]{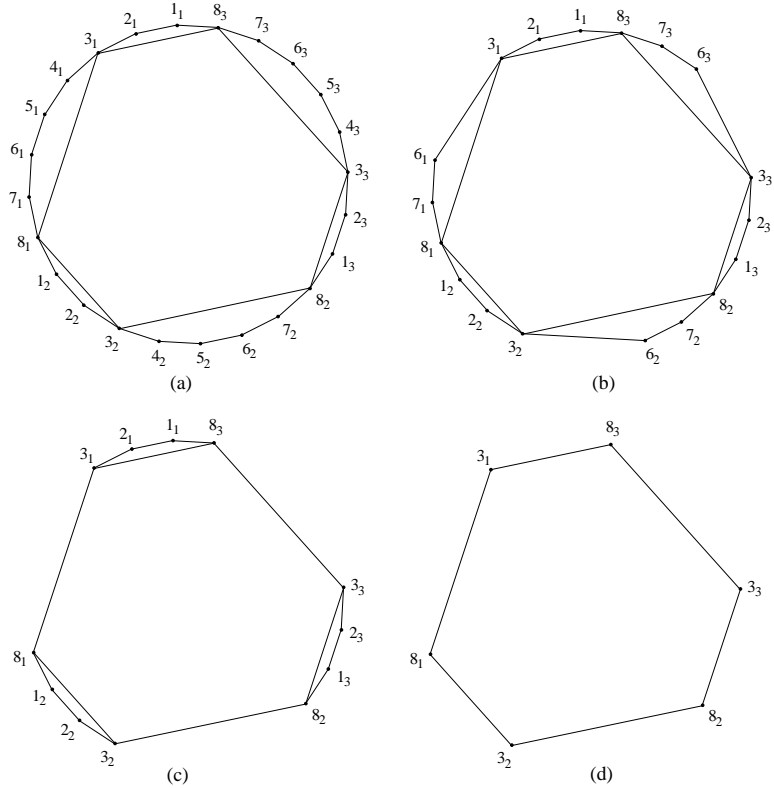}
\end{center}
\caption{\small A dissection of a 24-gon carried to the pair
$((3,8),(0,1,1))$ by Proposition \ref{pro:bijection}.}
\label{fig:dissection}
\end{figure}

\medskip
For illustration, Figure \ref{fig:dissection}(a) is a dissection of
a polygon with 24 vertices ($s=2$, $n=11$) using $6$ noncrossing
diagonals, which is invariant under 3-fold rotation $C_3$. These
diagonals are partitioned into two orbits under $C_3$ with labels 3
and 8, respectively. The bijection in Proposition
\ref{pro:bijection} carries the dissection to the pair of sequences
$(a_1,a_2)=(3,8)$ and $(\epsilon_1,\epsilon_2,\epsilon_3)=(0,1,1)$
as shown in (b)-(d) of Figure \ref{fig:dissection}.

\medskip
\begin{cor} \label{cor:U(s,n,k,d)} For $d\ge 2$ a common divisor of $sn+2$ and $k$, we have
\[ |U(s,n,k,d)|=
{\displaystyle
{{\frac{sn+2+k}{d}-1}\choose{\frac{k}{d}}}{{\lfloor\frac{n-1}{d}\rfloor}\choose{\frac{k}{d}}}}.
\]
\end{cor}

\medskip Since the results of Proposition
\ref{pro:k-odd} and Corollary \ref{cor:U(s,n,k,d)} agree with that
of Proposition \ref{pro:evaluation}, the proof of Theorem
\ref{thm:CSP} is completed.

\bigskip Of independent interest, we notice that the structure of
noncrossing trees appears to admit the cyclic sieving phenomenon,
which can be recovered immediately from the one for quadrangulations
of a $(2n+2)$-gon. A {\em noncrossing tree} is a tree drawn on a
circle of $n+1$ points numbered counterclockwise from 1 to $n+1$
such that the edges lie entirely within the circle and do not cross.
We shall establish a bijection $\Lambda$ between the set $X$ of
noncrossing trees with $n+1$ vertices and the set $Y$ of
quadrangulations of a $(2n+2)$-gon $P$, with the cyclic sieving
phenomenon preserved.

Given a $\pi\in X$, we associate the edge set $\{E_1,\dots,E_{n}\}$
of $\pi$ with a set of noncrossing lines $\{L_1,\dots,L_{n}\}$
within $P$ such that if $E_i$ connects points $j,k\in [n+1]$ of
$\pi$ then $L_i$ connects points $2j-1$, $2k-1$ of $P$, for each $i$
($1\le i\le n$). Then $\Lambda(\pi)$ is defined to be the unique
quadrangulation of $P$ that has $\{L_1,\dots,L_{n}\}$ as the set of
diagonals of all quadrilaterals. Conversely, given a $\pi'\in Y$, we
observe that each quadrilateral of $\pi'$ has a diagonal the
endpoints of which are labeled with odd numbers. Hence
$\Lambda^{-1}$ is obtained simply by a reverse procedure. Under this
bijection, rotation on a noncrossing tree $\pi$ by one point
clockwise corresponds to rotation on $\Lambda(\pi)$ by two points
clockwise.

\begin{figure}[ht]
\begin{center}
\includegraphics[width=3.6in]{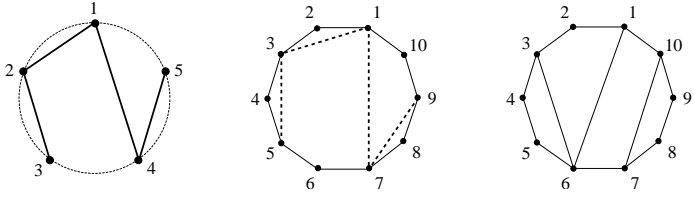}
\end{center}
\caption{\small A noncrossing tree on $5$ points, and the
corresponding quadrangulation of a $10$-gon.}
\label{fig:noncrossing}
\end{figure}

\medskip For example, on the left of Figure \ref{fig:noncrossing} is
a noncrossing tree with 5 points. The corresponding quadrangulation
of a $10$-gon is shown on the right,  the set of diagonals of all
quadrilaterals of which is shown in the center.

By Theorem \ref{thm:CSP} with $(s,k)=(2,n-1)$ and bijection
$\Lambda$, we have the following result.

\medskip
\begin{thm} Let $X$ be the set of noncrossing trees with $n+1$, and
let the cyclic group $C$ of order $n+1$ act on $X$ by cyclic
rotation of the vertices. Let
$X(q):=\frac{1}{[2n+1]_q}{{3n}\brack{n}}_q$. Then $(X,X(q),C)$
exhibits the cyclic sieving phenomenon.
\end{thm}

We remark that $X(1)$ recovers a known formula, due to M. Noy
\cite{N}, for the number of noncrossing trees with $n+1$ vertices.

\medskip
\section{The cyclic sieving phenomenon for $\Delta^s(B_n)$}

Following \cite[Section 5.2]{FR}, the generalized cluster complex
$\Delta^s(\Phi)$ for type $B_n$ can be realized as follows (see also
\cite{Simion, Tz-JCTA}). Let $P$ be a regular polygon with $2sn+2$
vertices. The vertices are labeled by
$\{1,2,\dots,sn+1,\overline{1},\overline{2},\dots,\overline{sn+1}\}$
counterclockwise. A $B$-{\em diagonal} of $P$ is either (i) a {\em
diameter}, i.e., a diagonal that connects a pair of antipodal points
$i,\overline{i}$, for some $1\le i\le sn+1$, or (ii) a pair of
$s$-divisible diagonals $ij$, $\overline{i} \overline{j}$ for two
distinct
$i,j\in\{1,2,\dots,sn+1,\overline{1},\overline{2},\dots,\overline{sn+1}\}$,
nonconsecutive around the boundary of the polygon. (It is understood
that if $a=\overline{i}$, then $\overline{a}=i$). Note that a
$B$-diagonal dissects $P$ into a pair of $(sm+2)$-gons and a
centrally symmetric $(2s(n-m)+2)$-gon ($1\le m\le n$). The vertices
of the complex $\Delta^s(B_n)$ correspond to the $B$-diagonals of
$P$, and the faces of $\Delta^s(B_n)$ correspond to $s$-divisible
dissections of $P$ using $B$-diagonals. The maximal faces correspond
to centrally symmetric dissections of $P$ into $(s+2)$-gons. For
$s=1$, this complex is the dual complex of the $n$-dimensional {\em
cyclohedron}, or {\em Bott-Taubes polytope} (see \cite[Lecture
3]{FR2}). Under this bijection, the map $\Gamma_s$ corresponds to
clockwise rotation of $P$ carrying point $2$ to point $1$, etc.

Taking Gaussian coefficients with base $q^2$ in Theorem
\ref{thm:formula}(ii),  we define
\begin{equation} \label{eqn:H(q)}
H(s,n,k;q)={{sn+k}\brack{k}}_{q^2}{{n}\brack{k}}_{q^2},
\end{equation}
for $0\le k\le n$. Note that  $H(s,n,k;1)=f_k(B_n,s)$.  We prove
that the faces of the generalized cluster complex $\Delta^s(B_n)$
exhibit the cyclic sieving phenomenon under the group action
$\langle\Gamma_s\rangle$.

\medskip
\begin{thm} \label{thm:CSP-B} For positive integers $s$ and $n$, let $X$ be the set of $s$-divisible
dissections of a $(2sn+2)$-gon using $k$ noncrossing $B$-diagonals.
Let the cyclic group $C$ of order $2sn+2$ act on $X$ by cyclic
rotation of the polygon. Let $X(q):=H(s,n,k;q)$. Then $(X,X(q),C)$
exhibits the cyclic sieving phenomenon.
\end{thm}

\begin{figure}[ht]
\begin{center}
\includegraphics[width=4.4in]{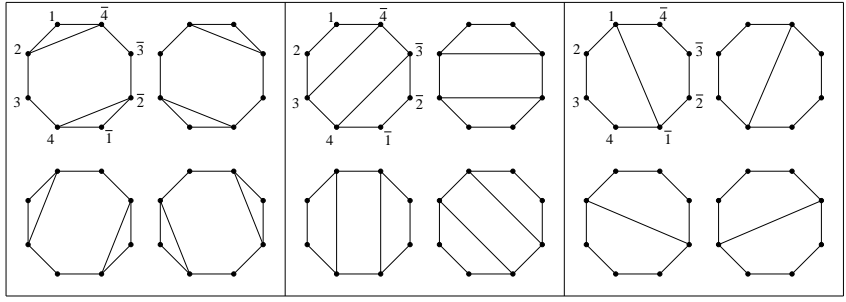}
\end{center}
\caption{\small The 1-divisible dissections of an octagon using 1
$B$-diagonal.} \label{fig:type-B}
\end{figure}

\medskip For example, take $s=1$, $n=3$, and $k=1$. $X(q)\equiv
3+3q^2+3q^4+3q^6$ (mod $q^8-1$). As shown in Figure
\ref{fig:type-B}, there are 12 1-divisible dissections of an octagon
using 1 $B$-diagonal. These dissections are partitioned into three
orbits under a group action by cyclic rotation, all of which have a
stabilizer of order 2.

For $t\ge 2$ a divisor of $sn+1$, let $\omega$ be a primitive $t$-th
root of unity. Making use of Lemma \ref{lem:Lemma} and Corollary
\ref{cor:d|k}, it is straightforward to prove that
\begin{equation} \label{eqn:B}
\lim_{q\rightarrow\omega} {{sn+k}\brack{k}}_q{{n}\brack{k}}_q
 = \left\{
\begin{array}{cl} \vspace{0.15cm}
                                {{\frac{sn+1+k}{t}-1}\choose{\frac{k}{t}}}{{\lfloor\frac{n-1}{t}\rfloor}\choose{\frac{k}{t}}} & \mbox{if
                                $t|k$,} \\
                                0  & \mbox{otherwise.}
           \end{array}
   \right.
\end{equation}
Note that
${\lfloor\frac{n-1}{t}\rfloor}={\lfloor\frac{n}{t}\rfloor}$ since
$t$ divides $sn+1$.

\medskip
\begin{pro} \label{pro:evaluation-B}
For $d\ge 2$ a divisor of $2sn+2$, let $\omega$ be a primitive
$d$-th root of unity. Then
\[
\lim_{q\rightarrow\omega}
{{sn+k}\brack{k}}_{q^2}{{n}\brack{k}}_{q^2}=\left\{
\begin{array}{cl}

\vspace{0.15cm}

{\displaystyle {{sn+k}\choose{k}}{{n}\choose{k}}}
& \mbox{if $d=2$,}\\

\vspace{0.15cm}
 {\displaystyle {{\frac{sn+1+k}{d}-1}\choose{\frac{k}{d}}}{{\lfloor\frac{n-1}{d}\rfloor}\choose{\frac{k}{d}}}} & \mbox{if
                          $d\ge 3$ odd, $d|k$,} \\
\vspace{0.15cm}
 {\displaystyle {{\frac{2sn+2+2k}{d}-1}\choose{\frac{2k}{d}}}{{\lfloor\frac{2(n-1)}{d}\rfloor}\choose{\frac{2k}{d}}}} & \mbox{if
                          $d\ge 4$ even, $d|2k$,} \\

0 & \mbox{otherwise.}
\end{array}
 \right.
\]
\end{pro}

\medskip
\begin{proof} For $d=2$, $\lim_{q\rightarrow -1}
{{sn+k}\brack{k}}_{q^2}{{n}\brack{k}}_{q^2}={{sn+k}\choose{k}}{{n}\choose{k}}$.
For $d\ge 3$, let $t=d$ if $d$ is odd, and $t=\frac{d}{2}$
otherwise. Then $t$ divides $sn+1$ and $\omega'=\omega^2$ is a
primitive $t$-th root of unity. By (\ref{eqn:B}), we have
\[\lim_{q\rightarrow\omega}
{{sn+k}\brack{k}}_{q^2}{{n}\brack{k}}_{q^2}=\lim_{q^2\rightarrow\omega'}
{{sn+k}\brack{k}}_{q^2}{{n}\brack{k}}_{q^2}= \left\{
\begin{array}{cl} \vspace{0.15cm}
                                {{\frac{sn+1+k}{t}-1}\choose{\frac{k}{t}}}{{\lfloor\frac{n-1}{t}\rfloor}\choose{\frac{k}{t}}} & \mbox{if
                                $t|k$,} \\
                                0  & \mbox{otherwise.}
\end{array}
\right.
\]
The assertion follows.
\end{proof}

\medskip
For $d\ge 2$ a divisor of $2sn+2$, let $C_d$ be a subgroup of order
$d$ of $C$ and let $V(s,n,k,d)\subseteq X$ denote the set of
dissections that are invariant under $C_d$. Note that $V(s,n,k,2)=X$
since all dissections in $X$ are invariant under 2-fold rotation.
For $d\ge 3$, we observe that if $V(s,n,k,d)$ is nonempty then each
$\pi\in V(s,n,k,d)$ consists of $k$ non-diameter $B$-diagonals, and
hence $k\le n-1$ and $d|2k$. Moreover, $V(s,n,k,d)=V(s,n,k,2d)$ if
$d$ is odd. The following result can be proved by the same argument
as the one in the proof of Proposition \ref{pro:bijection} and the
proof is omitted.

\medskip
\begin{pro} \label{pro:bijection-B} For $t\ge 2$ a common divisor of $sn+1$ and $k$,
there is a bijection between the set $V(s,n,k,2t)$ and the cartesian
product of the set of sequences $\{1\le a_1\le a_2\le\cdots\le
a_\frac{k}{t}\le\frac{sn+1}{t}\}$ and the set of sequences
$(\epsilon_1,\dots,\epsilon_m)\in\{0,1\}^m$ with exactly
$\frac{k}{t}$ entries equal to $1$, where
$m={\lfloor\frac{n-1}{t}\rfloor}$.
\end{pro}

\begin{figure}[ht]
\begin{center}
\includegraphics[width=4in]{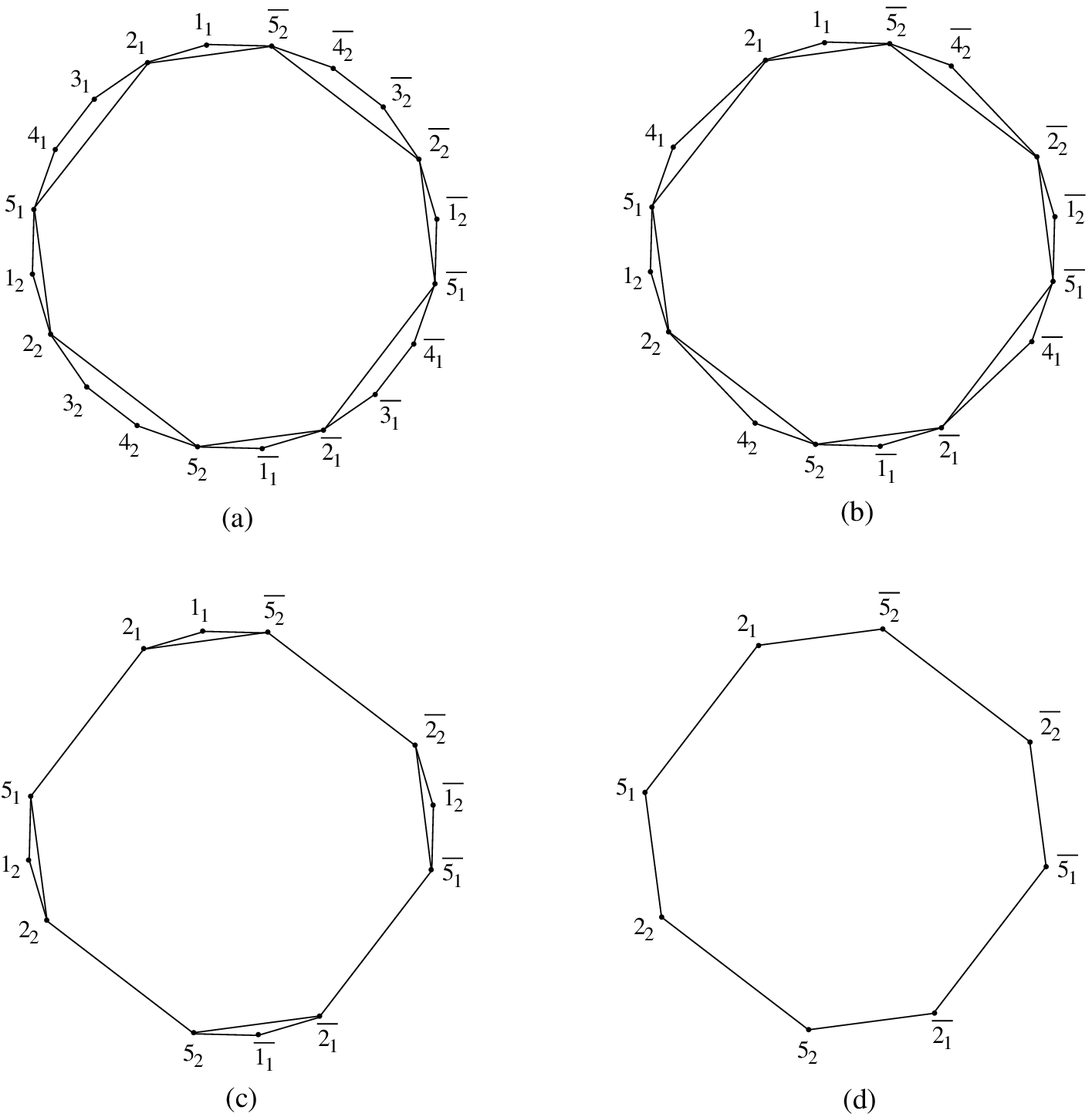}
\end{center}
\caption{\small A dissection of a 20-gon carried to the pair
$((2,5),(0,1,1,0))$ by Proposition \ref{pro:bijection-B}.}
\label{fig:B-dissection}
\end{figure}

\medskip
For illustration, Figure \ref{fig:B-dissection}(a) is a dissection
$\pi\in V(1,9,4,4)$ of a polygon with 20 vertices ($t=2$).  The
bijection in Proposition \ref{pro:bijection-B} carries $\pi$ to the
pair of sequences $(a_1,a_2)=(2,5)$ and
$(\epsilon_1,\epsilon_2,\epsilon_3,\epsilon_4)=(0,1,1,0)$ as shown
in (b)-(d) of Figure \ref{fig:B-dissection}.

\medskip Hence we have the following results.

\medskip
\begin{cor} \label{cor:V(s,n,k,d)} For $d\ge 2$ a divisor of $2sn+2$, we have
\[ |V(s,n,k,d)|=
\left\{ \begin{array}{cl}

\vspace{0.15cm}

{\displaystyle {{sn+k}\choose{k}}{{n}\choose{k}}} & \mbox{if $d=2$,}\\

\vspace{0.15cm}

{\displaystyle{{\frac{sn+1+k}{d}-1}\choose{\frac{k}{d}}}
                {{\lfloor\frac{n-1}{d}\rfloor}\choose{\frac{k}{d}}}} & \mbox{if $d\ge 3$ odd, $d|k$}\\

\vspace{0.15cm}

{\displaystyle{{\frac{2sn+2+2k}{d}-1}\choose{\frac{2k}{d}}}
                {{\lfloor\frac{2(n-1)}{d}\rfloor}\choose{\frac{2k}{d}}}} & \mbox{if $d\ge 4$ even, $d|2k$}\\
\vspace{0.15cm}

 0 & \mbox{otherwise,}
       \end{array}
\right.
\]
\end{cor}

Since the results of Corollary \ref{cor:V(s,n,k,d)} agree with that
of Proposition \ref{pro:evaluation-B}, the proof of Theorem
\ref{thm:CSP-B} is completed.

\section{The cyclic sieving phenomenon for $\Delta^s(D_n)$}

\begin{figure}[ht]
\begin{center}
\includegraphics[width=3.8in]{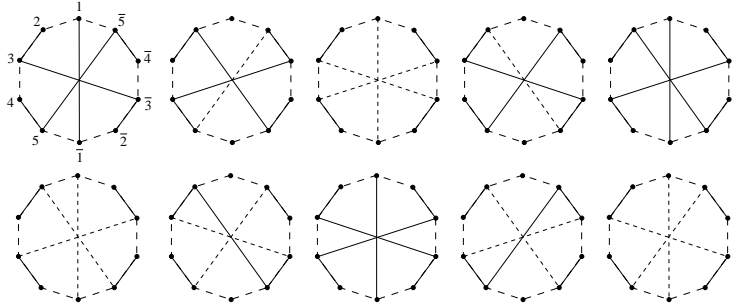}
\end{center}
\caption{\small The orbit of a maximal face of $\Delta^2(D_3)$.}
\label{fig:switcher}
\end{figure}

\medskip Following \cite[Section 5.3]{FR}, the generalized
cluster complex $\Delta^s(\Phi)$ for type $D_n$ can be realized as
follows. Let $P$ be a regular polygon with $2s(n-1)+2$ vertices. The
vertices are labeled by
$\{1,2,\dots,s(n-1)+1,\overline{1},\overline{2},\dots,\overline{s(n-1)+1}\}$
counterclockwise.  There are two copies of each diameter, one
colored red and the other colored blue. A vertex of the complex
$\Delta^s(D_n)$ is identified with a $D$-{\em diagonal} of $P$,
which is either (i) a red or a blue diameter, or (ii) a non-diameter
$B$-diagonal.  Let $\kappa(L_i)$ denote the color of $L_i$. The map
$\Gamma_s$ acts on $\Delta^s(D_n)$ by  rotating $P$ clockwise,
carrying point 2 to point 1, and switching the colors of certain
diameters. Specifically, $\Gamma_s$ carries the order pair
$(L_{j},\kappa(L_{j}))$ to $(L_{j-1},\kappa(L_{j-1}))$, where
$\kappa(L_{j-1})\neq\kappa(L_{j})$ if $j=1$ or $j\equiv 2$ (mod
$s$), and $\kappa(L_{j-1})=\kappa(L_{j})$ otherwise. The map
$\Gamma_s$ induces a relation of compatibility among $D$-diagonals.
Two diameters with the same endpoints and different colors are
compatible. Two diameters with distinct endpoints are compatible if
and only if applying $\Gamma_s$ repeatedly until either of them is
carried to $L_1$ results in diameters of the same color. In all the
other cases, two $D$-diagonals are compatible if they are
noncrossing in the sense of type-B dissections. For convenience, the
set $\{(sj+2,sj+1)|0\le j\le n-1\}$ of edges of $P$ are called {\em
color-switchers}, where $s(n-1)+2=\overline{1}$. Figure
\ref{fig:switcher} shows the orbit of a maximal face of
$\Delta^2(D_3)$ under the action of $\Gamma_2$, along with the
color-switchers, drawn as broken edges, indicating the locations at
which diameters change colors.

To describe the feasible polynomial $X(q)$ that serves our purpose,
we rewrite the $k$-face number of $\Delta^s(D_n)$ (Theorem
\ref{thm:formula}(iii)) as the sum of the expressions
\begin{equation}\label{eqn:D-1}
 {\displaystyle
{{s(n-1)+k}\choose{k}}{{n-1}\choose{k}}+{{s(n-1)+k}\choose{k}}{{n-2}\choose{k-1}}}
\end{equation}
and \begin{equation}\label{eqn:D-2}
 {\displaystyle
{{s(n-1)+k-1}\choose{k-1}}{{n-2}\choose{k-2}}+2{{s(n-1)+k-1}\choose{k}}{{n-2}\choose{k-2}}}.
\end{equation}
Note that the first term of (\ref{eqn:D-1}) is $f_k(B_{n-1},s)$,
which counts the $k$-faces of $\Delta^s(D_n)$ with at most one
diameter (say colored red). By a result of Tzanaki \cite[Corollary
3.4]{Tz-JCTA}, the second term of (\ref{eqn:D-1}) is the number of
$k$-faces of $\Delta^s(B_{n-1})$ with one diameter, which counts the
$k$-faces of $\Delta^s(D_n)$ with exactly one diameter (say colored
blue). Hence (\ref{eqn:D-1}) is the number of $k$-faces of
$\Delta^s(D_n)$ with at most one diameter of either color, and the
difference of the two terms is equal to the number of $k$-faces of
$\Delta^s(D_n)$ without diagonals. On the other hand, the first term
of (\ref{eqn:D-2}) is the number of $(k-1)$-faces of
$\Delta^s(B_{n-1})$ with one diagonal, which also counts the
$k$-faces of $\Delta^s(D_n)$ with two diagonals of different color
and the same endpoints. The second term of (\ref{eqn:D-2}) counts
the remaining $k$-faces of $\Delta^s(D_n)$, i.e., the ones with at
least two diagonals connecting distinct endpoints. Such faces can be
paired off by switching the colors of each diameters. For the sake
of being consistent with the polynomial $\Cat^{(s)}(D_n,q)$ in
(\ref{eqn:q-facet}) when we take $q$-analogues, (\ref{eqn:D-3}) is
alternatively written as
\begin{equation}\label{eqn:D-3}
 {\displaystyle
{{s(n-1)+k}\choose{k}}{{n-2}\choose{k-2}}+{{s(n-1)+k-1}\choose{k}}{{n-2}\choose{k-2}}}
\end{equation}
so that the difference of the two terms is equal to the number of
$k$-faces of $\Delta^s(D_n)$ with two diagonals of different color
and the same endpoints.

As a $q$-analog of the sum of (\ref{eqn:D-1}) and (\ref{eqn:D-3}),
we define the polynomial $F(s,n,k;q)$ by
\begin{equation}\label{eqn:D(q)}
 \begin{array}{l} \vspace{0.15cm}
F(s,n,k;q):={\displaystyle
{{s(n-1)+k}\brack{k}}_{q^2}{{n-1}\brack{k}}_{q^2}+{{s(n-1)+k}\brack{k}}_{q^2}{{n-2}\brack{k-1}}_{q^2}}\cdot
q^n
\\ \vspace{0.15cm} \hspace{3cm} {\displaystyle
+{{s(n-1)+k}\brack{k}}_{q^2}{{n-2}\brack{k-2}}_{q^2}+{{s(n-1)+k-1}\brack{k}}_{q^2}{{n-2}\brack{k-2}}_{q^2}}\cdot
q^n.
\end{array}
\end{equation}
Note that $F(s,n,k;1)=f_k(D_n,s)$. We shall prove that the $k$-faces
of the generalized cluster complex $\Delta^s(D_n)$, along with
$F(s,n,k;q)$, exhibit the cyclic sieving phenomenon.

\medskip
\begin{thm} \label{thm:CSP-D} For positive integers $s$ and $n$, let $X$ be the set of $s$-divisible
dissections of a $(2s(n-1)+2)$-gon using $k$ compatible
$D$-diagonals. Let $C$ be the cyclic group of order $2s(n-1)+2$
generated by $\Gamma_s$ that acts on $X$. Let $X(q):=F(s,n,k;q)$.
Then $(X,X(q),C)$ exhibits the cyclic sieving phenomenon.
\end{thm}

\begin{figure}[ht]
\begin{center}
\includegraphics[width=5.85in]{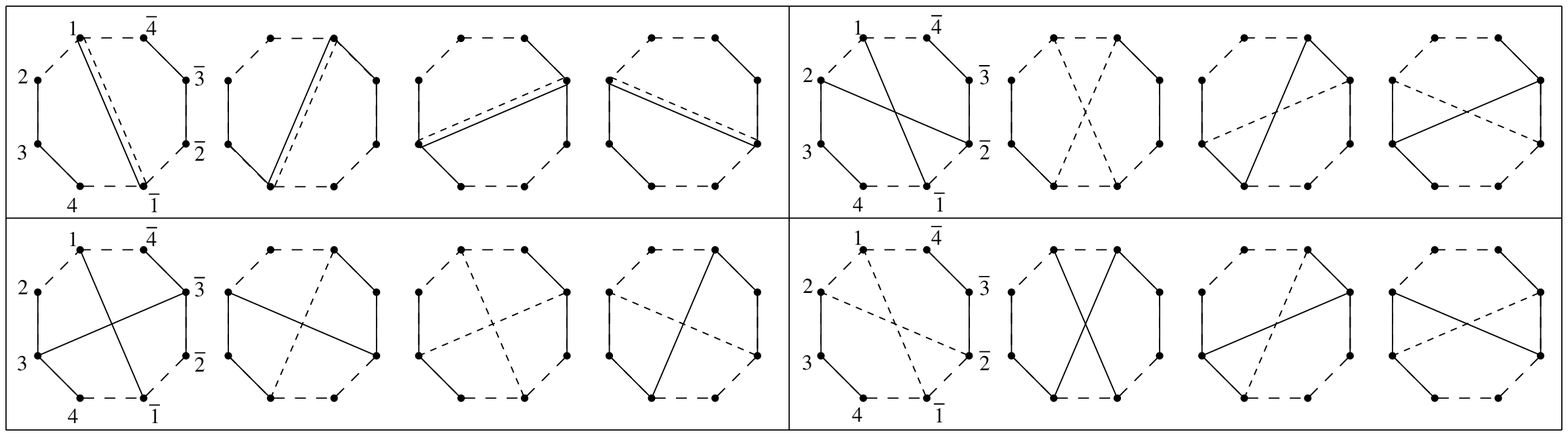}
\end{center}
\caption{\small The 3-divisible dissections of an octagon using 2
compatible $D$-diagonals.} \label{fig:type-D}
\end{figure}

\medskip For example, take $s=3$, $n=2$, and $k=2$. $X(q)\equiv 4+4q^2+4q^4+4q^6$
(mod $q^8-1$). As shown in Figure \ref{fig:type-D}, there are 16
$3$-divisible dissections of an octagon using 2 compatible
$D$-diagonals. These dissections are partitioned into four orbits
under the group action $\langle\Gamma_3\rangle$, all of which free.

To prove Theorem \ref{thm:CSP-D}, we first evaluate $F(s,n,k,q)$ at
a primitive $d$-th root $\omega$ of unity, where $d$ divides
$2s(n-1)+2$. For $d=2$, we have an immediate evaluation at $q=-1$.

\medskip
\begin{pro} \label{pro:d=2-in-X(q)}
\[ [F(s,n,k,q)]_{q=-1}=
\left\{ \begin{array}{cl} \vspace{0.15cm}

{\displaystyle
{{s(n-1)+k}\choose{k}}{{n}\choose{k}}+{{s(n-1)+k-1}\choose{k}}{{n-2}\choose{k-2}}}
& \mbox{if $n$ even,}\\
{\displaystyle
{{s(n-1)+k}\choose{k}}{{n-2}\choose{k}}+{{s(n-1)+k-1}\choose{k-1}}{{n-2}\choose{k-2}}}
& \mbox{otherwise.}
       \end{array}
\right.
\]
\end{pro}

For $d\ge 3$, by Lemma \ref{lem:Lemma} and Corollary \ref{cor:d|k},
we the have the following results as an intermediate stage of
evaluation.
\medskip
\begin{lem} \label{lem:intermedate} For $t\ge 2$ a divisor of $s(n-1)+1$, let $\omega$ be a  primitive
$t$-th root of unity. Then
\begin{enumerate}

\item \[\lim_{q\rightarrow\omega} {{s(n-1)+k}\brack{k}}_q{{n-1}\brack{k}}_q
=\left\{ \begin{array}{cl}
{{\frac{s(n-1)+1+k}{t}-1}\choose{\frac{k}{t}}}{{\lfloor\frac{n-2}{t}\rfloor}\choose{\frac{k}{t}}}
& \mbox{if $t|k$}\\
0 & \mbox{otherwise.}
\end{array}
 \right.
\]

\item \[\lim_{q\rightarrow\omega} {{s(n-1)+k}\brack{k}}_q{{n-2}\brack{k-1}}_q
=0
\]

\item \[\lim_{q\rightarrow\omega} {{s(n-1)+k}\brack{k}}_q{{n-2}\brack{k-2}}_q
=\left\{ \begin{array}{cl}
{{\frac{s(n-1)+1+k}{t}-1}\choose{\frac{k}{t}}}{{\frac{n}{t}-1}\choose{\frac{k}{t}-1}}
& \mbox{if $t|k$, $t|n$}\\
0 & \mbox{otherwise.}
\end{array}
 \right.
\]

\item \[
\lim_{q\rightarrow\omega}
{{s(n-1)+k-1}\brack{k}}_q{{n-2}\brack{k-2}}_q=\left\{\begin{array}{cl}
{{\frac{s(n-1)+1+k}{t}-1}\choose{\frac{k}{t}}}{{\frac{n}{t}-1}\choose{\frac{k}{t}-1}}
& \mbox{if $t|k$, $t|n$}\\
0 & \mbox{otherwise.}
\end{array}
 \right.
\]

\end{enumerate}

\end{lem}

\medskip
By Lemma \ref{lem:intermedate}, the following results can be proved
by an argument similar to the one in the proof of Proposition
\ref{pro:evaluation-B}.

\medskip
\begin{pro} \label{pro:evaluation-D}
For $d\ge 3$ a divisor of $2s(n-1)+2$, let $\omega$ be a  primitive
$d$-th root of unity. Then
\[ \begin{array}{l} \vspace{0.2cm}

[F(s,n,k;q)]_{q=\omega}= \\ \vspace{0.2cm}

\left\{ \begin{array}{cl}

\vspace{0.15cm} {\displaystyle
{{\frac{s(n-1)+1+k}{d}-1}\choose{\frac{k}{d}}}{{\frac{n}{d}}\choose{\frac{k}{d}}}
+{{\frac{s(n-1)+1+k}{d}-1}\choose{\frac{k}{d}}}{{\frac{n}{d}-1}\choose{\frac{k}{d}-1}}}
& \mbox{if $d\ge 3$ odd, $d|k$, $d|n$,}\\

\vspace{0.15cm} {\displaystyle
{{\frac{s(n-1)+1+k}{d}-1}\choose{\frac{k}{d}}}{{\lfloor\frac{n-2}{d}\rfloor}\choose{\frac{k}{d}}}}
& \mbox{if $d\ge 3$ odd, $d|k$, $d\nmid n$,}\\

\vspace{0.15cm} {\displaystyle
{{\frac{2s(n-1)+2+2k}{d}-1}\choose{\frac{2k}{d}}}{{\frac{2n}{d}}\choose{\frac{2k}{d}}}
+{{\frac{2s(n-1)+2+2k}{d}-1}\choose{\frac{2k}{d}}}{{\frac{2n}{d}-1}\choose{\frac{2k}{d}-1}}}
& \mbox{if $d\ge 4$ even, $d|2k$, $d|n$,}\\

\vspace{0.15cm} {\displaystyle
{{\frac{2s(n-1)+2+2k}{d}-1}\choose{\frac{2k}{d}}}{{\lfloor\frac{2(n-2)}{d}\rfloor}\choose{\frac{2k}{d}}}}
& \mbox{if $d\ge 4$ even, $d|2k$, $d\nmid n$,}\\

0 & \mbox{otherwise.}
\end{array}
 \right.

\end{array}
\]
\end{pro}

Let $X$ be the set of $s$-divisible dissections $\pi$ of $P$ using
$k$ compatible $D$-diagonals, for $0\le k\le n$. Note that if $\pi$
contains no diameter then $k\le n-2$. Let $C$ be the cyclic group of
order $2s(n-1)+2$ generated by $\Gamma_s$. For $d\ge 2$ a divisor of
$2s(n-1)+2$, let $W(s,n,k,d)\subseteq X$ denote the set of
dissections that are invariant under $d$-fold rotation $C_d\subseteq
C$. We observe that if $W(s,n,k,d)$ is nonempty then either $d=2$,
or $d\ge 3$ and $d|2k$. These two cases are treated in Propositions
\ref{pro:d=2-type-D} and
 \ref{pro:W(s,n,k,d)}, respectively.

\medskip
\begin{pro} \label{pro:d=2-type-D}
For positive integers $s$, $n$, and $k$,  we have
\[ |W(s,n,k,2)|=
\left\{ \begin{array}{cl} \vspace{0.15cm}

{\displaystyle
{{s(n-1)+k}\choose{k}}{{n}\choose{k}}+{{s(n-1)+k-1}\choose{k}}{{n-2}\choose{k-2}}}
& \mbox{if $n$ even,}\\
{\displaystyle
{{s(n-1)+k}\choose{k}}{{n-2}\choose{k}}+{{s(n-1)+k-1}\choose{k-1}}{{n-2}\choose{k-2}}}
& \mbox{otherwise.}
       \end{array}
\right.
\]
\end{pro}

\medskip
\begin{proof}
Note that, by the difference of the two terms of (\ref{eqn:D-1}),
there are ${{s(n-1)+k}\choose{k}}{{n-2}\choose{k}}$ dissections
without diameters, all of which are invariant under 2-fold rotation.
Given a dissection $\pi\in W(s,n,k,2)$ containing a diameter $L$,
when $P$ is rotated 180 degrees, $L$ passes $n$ color-switchers.
Then $L$ remains in the same color if $n$ is even, and switches to
the other color otherwise. Hence if $n$ is even, then any dissection
that contains diameters is invariant under 2-fold rotation, and
hence $|W(s,n,k,2)|=|X|$. Otherwise $n$ is odd and $\pi$ contains
two diameters of different colors with the same endpoints. By the
first term of (\ref{eqn:D-2}), there are
${{s(n-1)+k-1}\choose{k-1}}{{n-2}\choose{k-2}}$ such dissections,
and hence
$|W(s,n,k,2)|={{s(n-1)+k}\choose{k}}{{n-2}\choose{k}}+{{s(n-1)+k-1}\choose{k-1}}{{n-2}\choose{k-2}}$.
The assertion follows.
\end{proof}

For $d\ge 3$, we have the following necessary condition.

\medskip
\begin{lem} \label{lem:d|n} For $d\ge 3$ a common divisor of $2s(n-1)+2$
and $2k$, if there is a dissection in $W(s,n,k,d)$ containing a
diameter then $d$ divides $n$.
\end{lem}

\medskip
\begin{proof} Given a $\pi\in
W(s,n,k,d)$ containing a diameter $M$, let
$\M=\{M_1,\dots,\M_d\}\subseteq\pi$ be the orbit of diameters
containing $M$ under $C_d$. If $M_1=L_1$ then these diameters are of
the same color. We apply $\Gamma_s$ on $\pi$ repeatedly until
$M_{2}$ is carried to $M_{1}$ (i.e., $M_2=L_1$). Then each $M_i$
travels through $\frac{2s(n-1)+2}{d}$ edges. Let $M_i$ pass $r_i$
color-switchers. Clearly, $r_1+\cdots+r_d=2n$. Since color-switchers
are equally $s$ points apart, we have $|r_i-r_j|\le 1$, for $0\le
i,j\le d$. Moreover, the resulting diameters remain in the same
color, and this occurs whenever each $r_i$ is even. Consequently,
$r_1=\cdots=r_d=\frac{2n}{d}=2j$ for some $j$. Hence $d$ divides
$n$.

If $M_1\neq L_1$ then apply $\Gamma_s$ until $M_1=L_1$ so that the
diameters in $\M$ are of the same color, and then apply the argument
above. The proof is completed.
\end{proof}

\bigskip

We define the {\em initial point} of each $D$-diagonal. Each
diameter $L_i$ is oriented from $i$ to $\overline{i}$, for
$i\in[s(n-1)+1]$ and $i$ is called the initial point of $L_i$. Each
non-diameter $D$-diagonal $L$ is oriented in the usual sense, which
has two segments with starting points $j$, $\overline{j}$, for
$j\in[s(n-1)+1]$, and $j$ is called the initial point of $L$. The
maximal faces of $\Delta^s(D_n)$ consisting entirely of diameters
have the following properties, which are useful in enumerating
$W(s,n,k,d)$.

\medskip
\begin{lem} \label{lem:maximal-faces} {\rm (\cite[Lemma 9.1]{FR})} For $n\ge 3$,
given $n$ diameters of $P$, the following conditions are equivalent.
\begin{enumerate}
\item There exists a assignment of colors for each diameter such
that the colored diameters are pairwise compatible.

\item If $1\le a_1\le a_2\le\cdots\le a_n\le s(n-1)+1$ are the initial
points of the diameters, then $a_{j+1}-a_j\le s$, for $1\le j\le n$
(with the convention $a_{n+1}=a_1+s(n-1)+1)$.
\end{enumerate}
If these conditions hold, then there are exactly two ways to assign
colors to the diameters, which are related by switching the colors
of each of the $n$ diameters.
\end{lem}

\medskip
To enumerate $W(s,n,k,d)$, it is partitioned into three subsets
$T_0(s,n,k,d)$, $T_1(s,n,k,d)$, and $T_2(s,n,k,d)$, where
$T_0(s,n,k,d)$ is the subset of dissections that contain no
diameter, and $T_1(s,n,k,d)$ (resp. $T_2(s,n,k,d)$) is the subset of
dissections with diameters such that the first diameter is red
(resp. blue). Note that, by Lemma \ref{lem:d|n}, $T_1(s,n,k,d)$ and
$T_2(s,n,k,d)$ are nonempty if $d|n$ and that there is an immediate
bijection between $T_1(s,n,k,d)$ and $T_2(s,n,k,d)$ by switching the
colors of the diameters, thus $|T_1(s,n,k,d)|=|T_2(s,n,k,d)|$. In
the following, we make use of the same bijective method as the one
in the proof of Proposition \ref{pro:bijection-B} to enumerate
$T_1(s,n,k,d)$ and $T_0(s,n,k,d)$. With an analogous argument, a
proof is given for completeness.

\medskip
\begin{pro} \label{pro:bijection-D-ii} For $t\ge 2$ a common divisor of $s(n-1)+1$, $k$ and $n$,
there is a bijection between the set $T_1(s,n,k,2t)$ and the
cartesian product of the set of sequences $\{1\le a_1\le
a_2\le\cdots\le a_\frac{k}{t}\le\frac{s(n-1)+1}{t}\}$ and the set of
sequences $(\epsilon_1,\dots,\epsilon_m)\in\{0,1\}^m$ with exactly
$\frac{k}{t}$ entries equal to $1$, where $m={\frac{n}{t}}$ and
$\epsilon_m=1$.
\end{pro}

\begin{proof}
Fixing $t$, the vertices of the polygon are partitioned into $2t$
sectors of size $b=\frac{s(n-1)+1}{t}$ with a vertex-labeling
$\{1_1,2_1,\dots,b_1\},\dots,\{1_t,2_t,\dots,b_t\}$,
$\{\overline{1_1},\overline{2_1},\dots,\overline{b_1}\},\dots,\{\overline{1_t},\overline{2_t},\dots,\overline{b_t}\}$
in counterclockwise order. For each $\pi\in T_1(s,n,k,2t)$, the
$D$-diagonals of $\pi$ are partitioned into $\frac{k}{t}$ orbits of
size $t$ such that each orbit $R$ contains a unique $D$-diagonal $L$
the initial point of which lies in the first sector
$\{1_1,2_1,\dots,b_1\}$. We associate $R$ with the label $i\in[b]$
if the initial point of $L$ is $i_1$. Let
$\mu=(a_1,\dots,a_{\frac{k}{t}})$ be the sequence of labels
associated to the orbits of $D$-diagonals in $\pi$.

The bijection $\Psi_m$ is established by induction on
$m=\frac{n}{t}$. Recall that $1\le k\le n$ and $t$ divides $k$. If
$m=1$ then $\frac{k}{t}=1$, and $\Psi_1$ carries $\pi$ to the pair
$((a_1),(\epsilon_1))$, where $a_1$ is the label of the unique orbit
of diameters and $\epsilon_1=1$. Assume that $m\ge 2$. If
$a_{j+1}-a_j\le s$, for $1\le j\le \frac{k}{t}$ (with the convention
$a_{\frac{k}{t}+1}=a_1+b$) then by Lemma \ref{lem:maximal-faces} we
have $k=n$, so define $\Psi_{m}(\pi)=(\mu,\nu)$, where
$\nu=(1,\dots,1)$. Otherwise, locate the first $a_j$ such that
$a_j+s+1\le a_{j+1}$. We set $\epsilon_1=1$ if
$L=(a_j(a_j+s+1),\overline{a_j}(\overline{a_j}+s+1)$ is a
$D$-diagonal in $\pi$, and set $\epsilon_1=0$ otherwise. Let $R$ be
the $D$-diagonal-orbit containing $L$. The remaining
$\frac{k}{t}-\epsilon_1$ $D$-diagonal-orbits form a dissection of a
$(2s(n-t-1)+2)$-gon $P'$, obtained from $P$ by removing the
vertex-orbit of the set
$\{a_j+1,\dots,a_j+s,\overline{a_j}+1,\dots,\overline{a_j}+s\}$
under $C_{2t}$. By induction the remaining $\epsilon_i$'s are
determined.

To find $\Psi_m^{-1}$, given a pair of sequences
$\mu=(a_1,\dots,a_{\frac{k}{t}})$ and
$\nu=(\epsilon_1,\dots,\epsilon_m)$ with $\frac{k}{t}$ entries equal
to 1 and $\epsilon_m=1$, we retrieve the dissection
$\Phi_{m}^{-1}((\mu,\nu))$ as follows. If $m=1$ then $\frac{k}{t}=1$
and $\Psi_1^{-1}(\mu,\nu)$ consists of the unique orbit of diameters
with label $a_1$. Assume that $m\ge 2$. If $a_{j+1}-a_j\le s$, for
$1\le j\le \frac{k}{t}$, then by Lemma \ref{lem:maximal-faces} we
have $k=n$. Then $\Psi_m^{-1}(\mu,\nu)$ is defined to be the
dissection consisting of the diameter-orbits labeled by the $a_i$'s.
Otherwise, locate the first $a_j$ such that $a_j+s+1\le a_{j+1}$.
Then there is an orbit $R$ containing the $D$-diagonal
$L=(a_j(a_j+s+1),\overline{a_j}(\overline{a_j}+s+1))$ if
$\epsilon_1=1$ and there is not if $\epsilon_1=0$. Next, let
$\mu'=(a_1,\dots,a_{\frac{k}{t}})$ if $\epsilon_1=0$ and
$\mu'=(a_1,\dots,a_{j-1},a_{j+1},\dots,a_{\frac{k}{t}})$ otherwise,
and $\nu'=(\epsilon_2,\dots,\epsilon_{m})$. Use the pair
$(\mu',\nu')$ to inductively determine the set $\R$ of
$\frac{k}{t}-\epsilon_1$ $D$-diagonal-orbits in the polygon $P'$
described above. Then we define $\Phi_{m}^{-1}((\mu,\nu))$ to be the
dissection determined by $\R$ if $\epsilon_1=0$ and by $\R\cup\{R\}$
if $\epsilon_1=1$.
\end{proof}

\begin{figure}[ht]
\begin{center}
\includegraphics[width=4.4in]{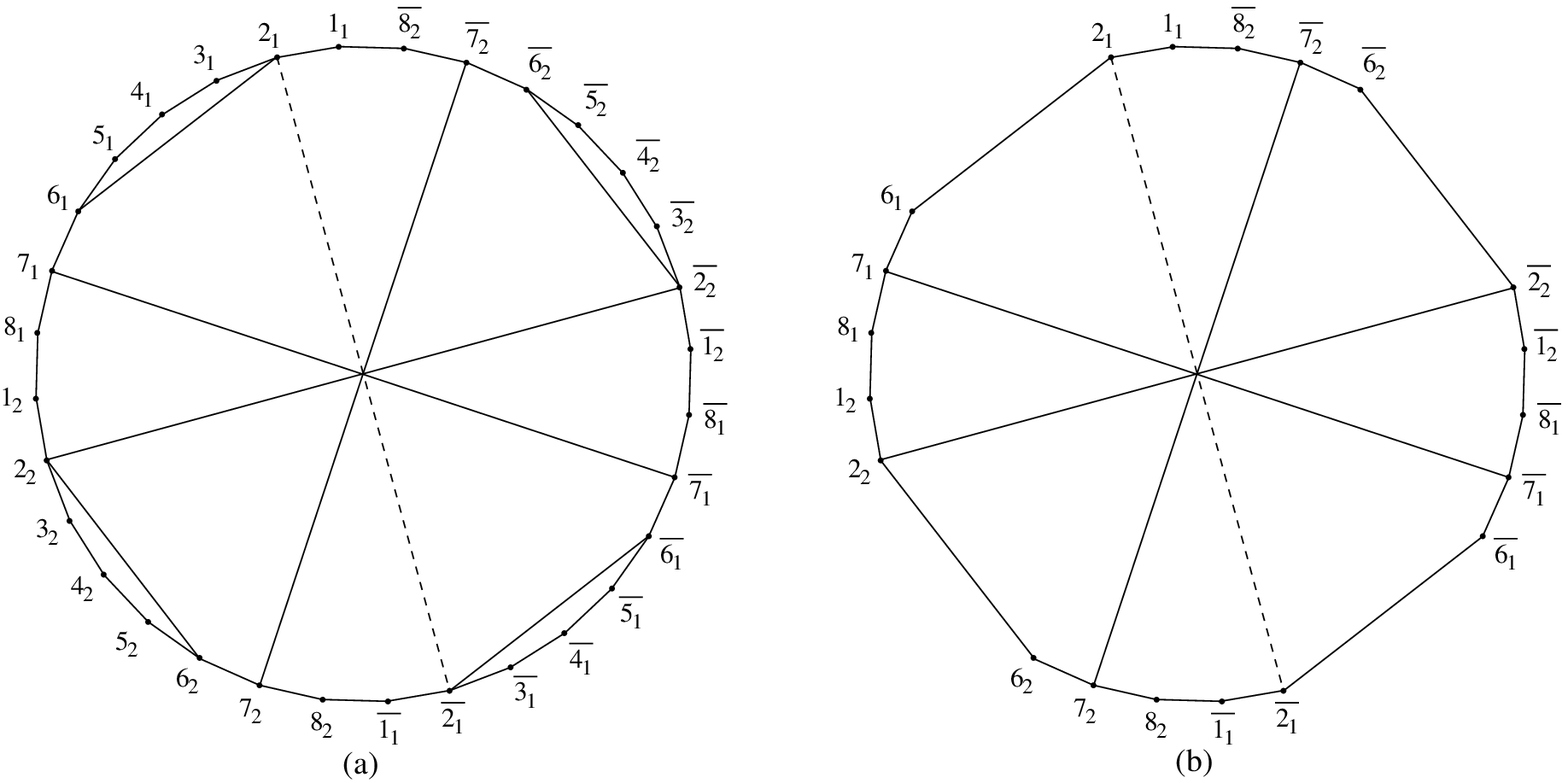}
\end{center}
\caption{\small A dissection of a 32-gon carried to the pair
$((2,2,7),(1,1,1))$ by Proposition \ref{pro:bijection-D-ii}}
\label{fig:D-dissection}
\end{figure}

For example, in Figure \ref{fig:D-dissection}(a), there is a
dissection $\pi\in T_1(3,6,6,4)$ of a polygon with 32 vertices
($t=2$). The bijection $\Psi_3$ in Proposition
\ref{pro:bijection-D-ii} carries $\pi$ to the pair of sequences
$(a_1,a_2,a_3)=(2,2,7)$ and
$(\epsilon_1,\epsilon_2,\epsilon_3)=(1,1,1)$ as shown in  (b) of
Figure \ref{fig:D-dissection}.

On the other hand, we observe that each dissection $\pi\in
T_0(s,n,k,2t)$ consists of $k$ non-diameter $D$-diagonals, and hence
$k\le n-2$ and $t|k$. As a case of Type-B dissection, the following
bijective result can be obtained by the same method as the one in
the proof of Proposition \ref{pro:bijection} and the proof is
omitted.

\begin{pro} \label{pro:bijection-D-i} For $t\ge 2$ a common divisor of $s(n-1)+1$ and $k$,
there is a bijection between the set $T_0(s,n,k,2t)$ and the
cartesian product of the set of sequences $\{1\le a_1\le
a_2\le\cdots\le a_\frac{k}{t}\le\frac{s(n-1)+1}{t}\}$ and the set of
sequences $(\epsilon_1,\dots,\epsilon_m)\in\{0,1\}^m$ with exactly
$\frac{k}{t}$ entries equal to $1$, where
$m={\lfloor\frac{n-2}{t}\rfloor}$.
\end{pro}

\medskip
By Propositions \ref{pro:bijection-D-ii} and
\ref{pro:bijection-D-i}, we have the following results.

\begin{pro} \label{pro:W(s,n,k,d)} For $d\ge 3$ a divisor of
$2s(n-1)+2$, we have
\[ \begin{array}{l} \vspace{0.2cm}

|W(s,n,k,d)|= \\ \qquad \vspace{0.2cm}

\left\{ \begin{array}{cl}

\vspace{0.15cm} {\displaystyle
{{\frac{s(n-1)+1+k}{d}-1}\choose{\frac{k}{d}}}{{\frac{n}{d}}\choose{\frac{k}{d}}}
+{{\frac{s(n-1)+1+k}{d}-1}\choose{\frac{k}{d}}}{{\frac{n}{d}-1}\choose{\frac{k}{d}-1}}}
& \mbox{if $d\ge 3$ odd, $d|k$, $d|n$,}\\

\vspace{0.15cm} {\displaystyle
{{\frac{s(n-1)+1+k}{d}-1}\choose{\frac{k}{d}}}{{\lfloor\frac{n-2}{d}\rfloor}\choose{\frac{k}{d}}}}
& \mbox{if $d\ge 3$ odd, $d|k$, $d\nmid n$,}\\

\vspace{0.15cm} {\displaystyle
{{\frac{2s(n-1)+2+2k}{d}-1}\choose{\frac{2k}{d}}}{{\frac{2n}{d}}\choose{\frac{2k}{d}}}
+{{\frac{2s(n-1)+2+2k}{d}-1}\choose{\frac{2k}{d}}}{{\frac{2n}{d}-1}\choose{\frac{2k}{d}-1}}}
& \mbox{if $d\ge 4$ even, $d|2k$, $d|n$,}\\

\vspace{0.15cm} {\displaystyle
{{\frac{2s(n-1)+2+2k}{d}-1}\choose{\frac{2k}{d}}}{{\lfloor\frac{2(n-2)}{d}\rfloor}\choose{\frac{2k}{d}}}}
& \mbox{if $d\ge 4$ even, $d|2k$, $d\nmid n$,}\\

0 & \mbox{otherwise.}
\end{array}
 \right.

\end{array}
\]
\end{pro}

\medskip
\begin{proof}
For $d\ge 3$, let $t=d$ if $d$ is odd and $t=\frac{d}{2}$ otherwise.
Note that $W(s,n,k,d)=W(s,n,k,2d)$ if $d$ is odd. If $t|k$ and $d|n$
then by Propositions \ref{pro:bijection-D-ii} and
\ref{pro:bijection-D-i}, we have
\begin{eqnarray*}
|W(s,n,k,2t)| &=& |T_0(s,n,k,2t)|+|T_1(s,n,k,2t)|+|T_2(s,n,k,2t)| \\
          &=& {{\frac{s(n-1)+1+k}{t}-1}\choose{\frac{k}{t}}}{{\frac{n}{t}-1}\choose{\frac{k}{t}}}
              +2{{\frac{s(n-1)+1+k}{t}-1}\choose{\frac{k}{t}}}{{\frac{n}{t}-1}\choose{\frac{k}{t}-1}}
              \\
          &=& {{\frac{s(n-1)+1+k}{t}-1}\choose{\frac{k}{t}}}{{\frac{n}{t}}\choose{\frac{k}{t}}}
+{{\frac{s(n-1)+1+k}{t}-1}\choose{\frac{k}{t}}}{{\frac{n}{t}-1}\choose{\frac{k}{t}-1}}.
\end{eqnarray*}
Moreover, if $t|k$ and $t\nmid n$ then by Lemma \ref{lem:d|n} and
Proposition \ref{pro:bijection-D-i}, we have
\[|W(s,n,k,2t)|=|T_0(s,n,k,t)|={{\frac{s(n-1)+1+k}{t}-1}\choose{\frac{k}{t}}}{{\lfloor\frac{n-2}{t}\rfloor}\choose{\frac{k}{t}}}. \]
Otherwise $t\nmid k$ and $|W(s,n,k,d)|=0$. The assertion follows.
\end{proof}

Since the results of Propositions \ref{pro:d=2-type-D} and
\ref{pro:W(s,n,k,d)} agree with that of Propositions
\ref{pro:d=2-in-X(q)} and \ref{pro:evaluation-D}, the proof of
Theorem \ref{thm:CSP-D} is completed.

\medskip
\noindent{\bf Remarks:} There is an alternative feasible option for
the polynomial $X(q)$ that is consistent with $\Cat^{(s)}(D_n,q)$
when $k=n$. The expression (\ref{eqn:D-1}) is rewritten as
\begin{equation}\label{eqn:D-4}
 {\displaystyle
{{s(n-1)+k}\choose{k}}{{n-2}\choose{k}}+2{{s(n-1)+k}\choose{k}}{{n-2}\choose{k-1}}},
\end{equation}
and then define $X(q)$, as a $q$-analogue of the sum of
(\ref{eqn:D-3}) and (\ref{eqn:D-4}), by
\begin{equation}\label{eqn:C(q)}
 \begin{array}{l} \vspace{0.15cm}
X(q):={\displaystyle
{{s(n-1)+k}\brack{k}}_{q^2}{{n-2}\brack{k}}_{q^2}+{{s(n-1)+k}\brack{k}}_{q^2}{{n-2}\brack{k-1}}_{q^2}}\cdot
(1+q^n)
\\ \vspace{0.15cm} \qquad\qquad {\displaystyle
+{{s(n-1)+k}\brack{k}}_{q^2}{{n-2}\brack{k-2}}_{q^2}
+{{s(n-1)+k-1}\brack{k}}_{q^2}{{n-2}\brack{k-2}}_{q^2}}\cdot q^n.
\end{array}
\end{equation}

\medskip
\section{The cyclic sieving phenomenon for $\Delta^s(I_2(a))$}

\begin{figure}[ht]
\begin{center}
\includegraphics[width=2.8in]{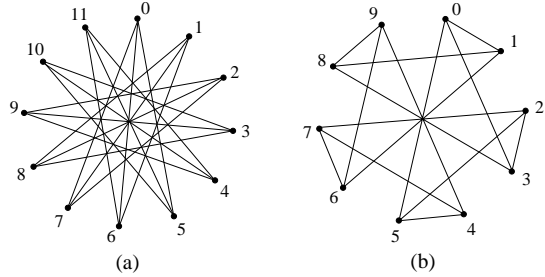}
\end{center}
\caption{\small Representation of $\Delta^2(I_2(5))$ and
$\Delta^2(I_2(4))$.} \label{fig:I2a}
\end{figure}

For a root system $\Phi$ of type $I_2(a)$, the complex
$\Delta^s(I_2(a))$ is an $(s+1)$-regular graph on $sa+2$ vertices,
the edges of which are facets. As shown in \cite[Example 4.4]{FR},
this graph can be constructed in the plane on a circle of $sa+2$
points labeled from 0 to $sa+1$ clockwise. For $a$ odd, the edge set
has $(am+2)$-fold rotational symmetry and connects each vertex $v$
to the $s+1$ vertices $v+\frac{s(a-1)}{2}+j$ (mod $sa+2$), for
$j=1,\dots,s+1$. Figure \ref{fig:I2a}(a) shows this graph for $s=2$
and $a=5$. In this case, the map $\Gamma_s$ corresponds to a
counterclockwise rotation of the graph by $\frac{2\pi}{sa+2}$. For
$a$ even, fixing an odd integer $i$, the edge set has
$(\frac{sa+2}{2})$-fold rotational symmetry and connects $0$ to the
vertices $i,i+2,\dots,i+2s$. Figure \ref{fig:I2a}(b) is
$\Delta^2(I_2(4))$ drawn in this style with $i=1$. In this case, the
map $\Gamma_s$ acts by a counterclockwise rotation of the graph by
$\frac{4\pi}{sa+2}$. Along with the $q$-analogue $X(q)$ of Theorem
\ref{thm:formula}(iv), we prove the cyclic sieving phenomenon for
the facets of $\Delta^s(\Phi)$.

\medskip
\begin{pro} \label{thm:CSP-I2(a)}
Let $X$ be the edge set of the graph $\Delta^s(I_2(a))$. Define
\begin{equation}\label{eqn:I2a} X(q):=\frac{[sa+2]_q}{[2]_q}\cdot\frac{[sa+a]_q}{[a]_q}.
\end{equation} Let $C=\ZZ_{sa+2}$ be the cyclic group that acts on $X$ by
cyclic rotation of the graph. Then $(X,X(q),C)$ exhibits the cyclic
sieving phenomenon.
\end{pro}

\medskip
\begin{proof} For $d\ge 2$ a divisor of $sa+2$, let $\omega$ be a primitive $d$-th root of
unity and let $Y(s,a,d)\subseteq X$ be the set of edges that are
invariant under a subgroup of order $d$ of $C$. There are two cases.

For $a$ odd,  it is clear that $[X(q)]_{q=\omega}=0$ for all $d\ge
3$. If $d=2$ then $s$ is even and $[X(q)]_{q=-1}=\frac{sa+2}{2}$.

On the other hand, an edge of the graph $\Delta^s(I_2(a))$ is
invariant under $d$-fold rotation only if $d=2$. By the construction
of the graph, $Y(s,a,2)$ is nonempty if and only if $s$ is even, in
which case $Y(s,a,2)$ consists of the $\frac{sa+2}{2}$ edges
connecting antipodal pairs of vertices. This proves that
$(X,X(q),C)$ exhibits the cyclic sieving phenomenon for $a$ odd.

For $a$ even, it is clear that $[X(q)]_{q=\omega}=0$ for all $d\ge
3$ and $[X(q)]_{q=-1}=\frac{(sa+2)(s+1)}{2}$. On the other hand, we
remark that under the group action $C$, odd (resp. even) vertices of
$\Delta^s(I_2(a))$ are carried to odd (resp. even) vertices, and
each vertex travels two cycles. Since each edge connects an odd
vertex and an even vertex, an edge is invariant if and only if it is
rotated by a full cycle. Hence $Y(s,a,2)=X$ and $Y(s,a,d)=\emptyset$
for $d\ge 3$. This proves that $(X,X(q),C)$ exhibits the cyclic
sieving phenomenon for $a$ even.
\end{proof}

It is straightforward to prove that the vertex set of
$\Delta^s(I_2(a))$, along with the polynomial $X(q)=[sa+2]_q$ if $a$
is odd and $X(q)=[sa+2]_{q^2}$ otherwise, exhibits the cyclic
sieving phenomenon since all vertices form an orbit if $a$ is odd
and two orbits otherwise under $C$.

\section{The cases of exceptional types}

In this section, we consider the generalized cluster complexes
$\Delta^s(\Phi)$ of exceptional type  $E_6$, $E_7$, $E_8$, $F_4$,
$H_3$, and $H_4$.

When $k=n$, the polynomial $X(q)=\Cat(\Phi,q)$ define in
(\ref{eqn:q-facet1}) is a feasible $q$-analogue of the number of
facets of $\Delta(\Phi)$. In particular, we verify Theorem
\ref{thm:allCSP-facets} affirmatively for the complex $\Delta(\Phi)$
of exceptional type. To see this, with the Coxeter numbers and
exponents of $\Phi$ listed in Figure \ref{fig:table}, the polynomial
$X(q)$ is expanded as follows.
\begin{enumerate}

\item In type $E_6$, $X(q)\equiv
67+52q+67q^2+52q^3+\cdots+67q^{12}+52q^{13}$ (mod $q^{14}-1$).

\item In type $E_7$, $X(q)\equiv 416+416q^2+416q^4+\cdots+416q^{18}$ (mod
$q^{20}-1$).

\item In type $E_8$, $X(q)\equiv
1574+1562q^2+1572q^4+1562q^6+\cdots+1574q^{24}+1562q^{26}+1572q^{28}+1562q^{30}$
(mod $q^{32}-1$).

\item In type $F_4$, $X(q)\equiv 15+15q^2+15q^4+\cdots+15q^{12}$ (mod
$q^{14}-1$).

\item In type $H_3$, $X(q)\equiv 6+5q^2+5q^4+6q^6+5q^8+5q^{10}$ (mod
$q^{12}-1$).

\item In type $H_4$, $X(q)\equiv
18+17q^2+18q^4+17q^6+\cdots+18q^{28}+17q^{30}$ (mod $q^{32}-1$).

\end{enumerate}

As searched by a computer, the orbit-structures for the $k$-faces of
the cluster complex $\Delta(\Phi)$  under the cyclic group $C$
generated by $\Gamma$ are shown in Figure
\ref{fig:orbit-structures}. We write $a_1(b_1)$, $a_2(b_2),\dots$,
$a_t(b_t)$ for the orbit-structure of the $k$-faces that are
partitioned into $b_i$ orbits of size $a_i$, for $1\le i\le t$,  in
which case the number of $k$-faces is equal to
$a_1b_1+\cdots+a_tb_t$.

\begin{figure}[ht]
\begin{tabular}{|c|c|c|c|c|c|c|}
  \hline
  $k$  & $E_6$ & $E_7$ & $E_8$ & $F_4$ & $H_3$ & $H_4$ \\
  \hline \hline
  1    & 14(2), 7(2)    & 10(7)         & 16(8)          & 7(4)        & 6(3)   & 16(4)\\
   \hline
  2    & 14(26), 7(5)   & 10(97), 5(1)  & 16(149), 8(3)  & 7(19)       & 6(8)   & 16(21), 8(1)\\
   \hline
  3    & 14(104), 7(13) & 10(518)       & 16(1121)       & 7(30)       & 6(5), 2(1) & 16(35)\\
   \hline
  4    & 14(195), 7(18) & 10(1410), 5(1)& 16(4211), 8(3), 4(2)  & 7(15)&            & 16(17), 8(1)\\
   \hline
  5    & 14(171), 7(15) & 10(2020), 2(1)& 16(8778)       &             &            & \\
   \hline
  6    & 14(52), 7(15)  & 10(1456)      & 16(10230), 8(22)      &      &            & \\
   \hline
  7    &                & 10(416)       & 16(6270)       &             &            & \\
   \hline
  8    &                &               & 16(1562), 8(10), 4(2) &      &            & \\
  \hline
\end{tabular}
\caption{\small The orbit-structures of the $k$-faces of
$\Delta(\Phi)$ of exceptional types.} \label{fig:orbit-structures}
\end{figure}

We observe that the orbit-structures for the maximal faces of
$\Delta(\Phi)$ shown in Figure \ref{fig:orbit-structures} agree with
the expansions (i)-(vi) of $X(q)$ (mod $q^{h+2}-1$). Hence the
cyclic sieving phenomenon holds. Consequently, together with $s=1$ in the
Theorems \ref{thm:CSP}, \ref{thm:CSP-B}, \ref{thm:CSP-D}, and
\ref{thm:CSP-I2(a)} for $k=n$, Theorem
\ref{thm:allCSP-facets} is proved.

\section{Discussion and remarks}
In this paper, we study the cyclic sieving phenomenon for the faces
of the generalized cluster complexes $\Delta^s(\Phi)$ associated a
root systems $\Phi$. As readers can see, the $q$-analogue $X(q)$ of face numbers plays
an essential role.

First we note that when $k=n$, our polynomials $X(q)$ for the facets
of $\Delta^s(\Phi)$ agree with the generalized $q$-Catalan numbers
defined by
\begin{equation} \label{eqn:q-facet} \Cat^{(s)}(\Phi,q):=\prod_{i=1}^n
\frac{[sh+e_i+1]_q}{[e_i+1]_q},
\end{equation}
i.e., $\Cat^{(s)}(A_{n-1},q)=G(s,n,n)$,
$\Cat^{(s)}(B_n,q)=H(s,n,n)$, $\Cat^{(s)}(D_n,q)=F(s,n,n)$, and
$\Cat^{(s)}(I_2(a),q)$ is (\ref{eqn:I2a}). As a result, we prove the
following conjecture mentioned by Reiner-Stanton-White, for
$\Delta^s(\Phi)$ in types $A_n$, $B_n$, $D_n$, and $I_2(a)$.

\medskip
\begin{con} {\rm (\cite{RSW2})}\label{con:CSP-facets} For a positive integer $s$, let $X$ be the set of
facets of the generalized cluster complex $\Delta^s(\Phi)$. Let
$X(q)=\Cat^{(s)}(\Phi,q)$ be defined in (\ref{eqn:q-facet}). Let the
cyclic group $C$ of order $sh+2$ generated by $\Gamma_s$ act on $X$.
Then $(X,X(q),C)$ exhibits the cyclic sieving phenomenon.
\end{con}

\medskip
The cases other than facets are not determined for lack of feasible
polynomials $X(q)$, although we believe that the cyclic sieving
phenomenon still hold once $X(q)$ is found. In fact, we would like
to point out that the {\em genuine} $q$-analogue $X(q)$ of the
$k$-face numbers of the complex $\Delta^s(\Phi)$ is quite elusive,
as discussed below.

As readers have seen, what we have  (i.e., (\ref{eqn:G(q)}),
(\ref{eqn:H(q)}), and (\ref{eqn:D(q)})) are feasible polynomials
$X(q)$ that serve the purpose of the cyclic sieving phenomenon. For
the type $A_n$ and $B_n$, we are more or less confident that the
polynomial (\ref{eqn:G(q)}), (\ref{eqn:H(q)}) are the genuine
$q$-analogues of the number of facets of $\Delta^s(\Phi)$.

However, our $q$-analogue $X(q)$ in type $D_n$ (\ref{eqn:D(q)}) is
still not so satisfying in the sense that it seems to be
artificially tailored to serve the purpose of cyclic sieving
phenomenon. We remark that in type $D_n$ same obvious $q$-analogues
of the expression in Theorem \ref{thm:formula}(iii) is not a
feasible option. In fact, besides (\ref{eqn:D(q)}) we have found
several other options of $X(q)$ for type $D_n$, each of them
expressed in sums of several terms and serving cyclic sieving
phenomenon well.

A unified formula of $k$-face numbers of $\Delta^s(\Phi)$ is found
by Fomin and Reading \cite{FR} to be
\begin{equation} \label{eq:FRcf}
f_k(\Phi,s)=c(\Phi,k,s){n\choose k}\prod_{L(e_i)\le k}\frac{sh+e_i+1}{e_i+1},
\end{equation}
where $c(\Phi,k,s){n\choose k}$ and $L(e_i)$ depends on $\Phi$, see
\cite{FR} for more information. However, except for the results
presented in the previous sections, we can hardly derive any other
feasible $q$-analogue (for cyclic sieving phenomenon) from obvious
$q$-analogues of this formula, even for $s=1$.

For example, for $k=1$ and $s=1$ in type $E_6$, by \cite{FR} we have
$f_1(\Phi,1)={{n}\choose{1}}\frac{h+e_1+1}{e_1+1}$. Taking a natural
$q$-analogue, we have
$X(q)={{6}\brack{1}}_q\frac{[14]_q}{[2]_q}\equiv
3+3q+3q^2+\cdots+3q^{13}$ (mod $q^{14}-1$), which does not agree
with the orbit-structure (i.e., with two orbits of size 14 and two
orbits of size 7) of the vertex set of $\Delta(E_6)$ under $C$. What
worse is that sometimes the natural $q$-analogue is not even a
polynomial (e.g., taking $s=1$ and $k=5$ in type $E_6$).

In light of these, we are very interested in finding a genuine polynomial
 $X(q)$ that involves the Coxeter number and exponents of $\Phi$
in the sense of (\ref{eqn:q-facet}). Such a polynomial $X(q)$ should
not only lead to a unified result for the $k$-faces of the complex
$\Delta^s(\Phi)$ as the one in Conjecture \ref{con:CSP-facets} but
also be consistent with other combinatorial structures (e.g.
noncrossing partitions) in connection with Coxeter groups. We leave
it as an open problem.

\noindent{\bf Open Problem.} What is the genuine $q$-analogue $X(q)$ of the
$k$-face numbers of the complex $\Delta^s(\Phi)$?

\bigskip \bigskip \centerline{\sc Acknowledgments}

The authors would like to thank Victor Reiner and Dennis Stanton for
suggesting this problem and for very helpful comments. This paper is
written during the first author's (S.-P. Eu) visit to Department of
Mathematics, University of Minnesota  and the second author's (T.-S.
Fu) visit to DIMACS center, Rutgers University. Both authors thank
the institutes for their hospitality. The authors also thank John
Stembridge for making available his {\tt coxeter} package for {\tt
Maple}, which is helpful in obtaining the results in Figure
\ref{fig:orbit-structures}. The first author is partially supported
by TJ \& MY Foundation.

%%%%%%%%%%% bibliography

\end{document}